\documentclass[a4paper,11pt]{amsart} 

\usepackage{amsfonts}
\usepackage{amsthm}
\usepackage{amssymb}
\usepackage[leqno]{amsmath}
\usepackage[english]{babel} 
\usepackage[leqno]{amsmath}
\usepackage{amssymb,amsthm}
\usepackage{enumerate}

\newtheorem{thm}{Theorem}[section]
\newtheorem{lemma}[thm]{Lemma}
\newtheorem{prop}[thm]{Proposition}
\newtheorem{cor}[thm]{Corollary}

\newtheorem{cordef}[thm]{Corollary - Definition}
\theoremstyle{definition}
\newtheorem{defi}[thm]{Definition}
\newtheorem{nota}[thm]{}

\theoremstyle{remark}
\newtheorem{remark}[thm]{Remark}

\newtheorem{example}[thm]{Example}

\newcommand{\la}{\longrightarrow}
\newcommand{\ha}{\hookrightarrow}
\newcommand{\da}{\dashrightarrow}

\newcommand{\ov}{\overline}

\def\X{\mathcal X}
\def\Y{\mathcal Y}
\def\ZZ{\mathcal Z}
\def\L{\mathcal L}
\def\N{\mathcal N}
\def\NL{\overline{\mathcal N}}
\def\O{\mathcal O}
\def\Poi{\mathcal P}

\def\ner{{\operatorname{N}}(\jac)}
\def\nerd{{\operatorname{N}}(\jacd)}
\def\NN{\operatorname{N}}
\def\NX{{\operatorname{N}} _X}
\def\NXb{\overline{\NX}}
\def\NS{{\operatorname{N}} _{\XS}}
\def\NY{{\operatorname{N}} _Y}

\def\dcg{\Delta _X}
\def\dcY{\Delta _{\YS}^d}
\def\dcYY{\Delta _{\YS}^{d,1}}
\def\dcYd{\Delta _Y^d}
\def\dcYS{\Delta _{Y,S}^d}
\def\glu{{\bf {G}} _X}

\def\dcS{\Delta _{\XS}^{d-s}}

\def\LX{\Lambda _X}

\newcommand{\G}{\Gamma _X}
\def\comp{\mu }
\def\int{M_X}
\def\vectn{n_1,\ldots, n_{\gamma}}
\def\vectd{d_1,\ldots, d_{\gamma}}
\def\dom{{\operatorname B}}
\def\domX{{\operatorname B}_X}

\newcommand{\dX}{\dot{\X}}
\def\mc{\underline{c}_i}

\def\md{\underline{d}}
\def\mb{\underline{b}}

\newcommand{\g}{\Gamma}

\newcommand{\pr}[1]{\mathbb{P}^{#1}}

\newcommand{\Z}{\mathbb{Z}}
\newcommand{\Q}{\mathbb{Q}}

\newcommand{\Aut}{\operatorname{Aut}}

\newcommand{\Stab}{\operatorname{Stab}}

\newcommand{\Hilb}{\operatorname{Hilb}^{dt-g+1}_{\pr{r}}}

\newcommand{\Spec}{\operatorname{Spec}}

\newcommand{\im}{\operatorname{Im}}

\newcommand{\mdeg}{\operatorname{{\underline{de}g}}}

\newcommand{\mgb}{\ov{M}_g}
\newcommand{\cgbar}{{\mathcal C}_g}
\newcommand\poin{\overline{\L }}

\def\pdg{\overline{P}_{d,\,g}}
\def\pd{P^d_g}

\def\taut{\mathcal{D}}

\def\pdgst{\overline{\mathcal{P}}_{d,g}}
\def\pdst{{\mathcal{P}_{d,g}}}
\def\mgbst{\overline{\mathcal{M}}_g}

\def\pdfun{{\mathcal{P}}^d_g}
\def\pX{P^d_X}
\def\pXb{\overline{P^d_X}}
\def\pf{P^d_f}
\def\pffun{{\mathcal{P}}^d_f}

\def\pfb{\overline{P^d_f}}
\def\Picfun{\mathcal{P}ic}

\def\jacf{\Pic^{\underline 0}_f}
\def\jac{\Pic _K^0}

\def\jacd{\Pic _K^d}

\def\jacdh{\Pic _{K_1}^d}
\newcommand{\Pic}{\operatorname{Pic}}

\newcommand{\pfun}[1]{{\Picfun}_f^{#1}}
\newcommand{\picf}[1]{\Pic_f^{#1}}
\newcommand{\pich}[1]{\Pic_h^{#1}}
\newcommand{\picX}[1]{\Pic^{#1}X}

\newcommand{\twX}{\operatorname{Tw}\negthickspace_fX}

\newcommand{\XS}{X_S^{\nu}}
\newcommand{\YS}{Y_S}

\newcommand{\gen}{\X_K}

\begin{document}


\begin{center}
{\bf\large
N\'eron models and compactified Picard schemes \\
over the moduli stack of stable curves} \\
\bigskip

Lucia  Caporaso

\bigskip

\end{center}

\noindent{\it Abstract.} We construct  modular Deligne-Mumford stacks  $\pdst$ representable
over
$\mgbst$ parametrizing N\'eron models of Jacobians as follows. 
Let  $B$ be a smooth  curve and $K$ its function field,
 let
$\gen$ be a smooth genus-$g$ curve over $K$ admitting stable minimal model
over $B$. The N\'eron model $\NN(\Pic^d\gen)\to B$ 
is then the base change of $\pdst$ via the moduli map $B\la \mgbst$
of $f$, i.e.:  $\NN(\Pic^d\gen)\cong \pdst\times _{\mgbst}B$.
Moreover  $\pdst$ is  compactified by a Deligne-Mumford stack  over $\mgbst$,
giving a completion of N\'eron models naturally stratified in terms of N\'eron models.
{\small \tableofcontents}

\section{Introduction}

\begin{nota}{\it Problems and results}
The first goal of this paper is  a parametrization result for 
N\'eron models of  Jacobians of stable curves (Theorem~\ref{picner}).
A technical part of the argument
that yields results of independent interest  is the strengthening
of a  construction
of the compactified  Picard variety over $\mgbst$. 
A further outcome is  a geometrically meaningful
compactification of such N\'eron models. 
We proceed to discuss
all of that more precisely.

Let $K=k(B)$ be the field of rational functions
of a  nonsingular one-dimensional  scheme $B$  defined over an algebraically closed
field $k$. Let $\gen$ be a nonsingular connected projective curve 
of genus $g\geq 2$ over $K$ whose regular minimal model 
over $B$ is a family $f:\X \la B$ of stable curves.

For any integer $d$, denote by $\jacd:=\Pic^d
\gen$  the degree-$d$ Picard variety of $\gen$
(parametrizing line bundles of degree $d$ on $\gen$),  and let
$\nerd$ be its N\'eron model over $B$.
It is well known that  (since the total space $\X$ is nonsingular)
the fibers of $\nerd\la B$ over
the closed points of $B$ depend only on the corresponding fibers of $f$.

It makes therefore sense to ask the following question: does there exist a space over 
$\mgbst$, such that, for every $K$ and $\gen$ as above, $\nerd$ is 
the base change of such a space via the moduli map $B\la \mgbst$ associated to the
family $f$?

In this paper 
we give a positive answer to this question  for every 
 $g\geq 3$ 
and for every $d$ such that $(d-g+1,2g-2)=1$.
Let us first state a result in scheme theoretic terms,
postponing the  stack-theoretic generalization for a moment
(cf. Theorem~\ref{picner}).
We construct
a  separated scheme $\pd$ over the moduli scheme of stable curves $\mgb$,
having the following property: for any
family $f:\X
\la
B$ of automorphism-free stable curves 
with  $\X$  regular, there is a canonical isomorphism of $B$-schemes
$$
\nerd \cong_B B\times_{\mgb}\pd
$$
where $B$ is viewed as an $\mgb$- scheme via the moduli map of the family $f$.

Working within the
category
of schemes,
the  restriction to automorphism-free curves is  necessary: if
$X$ is a stable curve, $\Aut (X)$ injects into the automorphism group of its
generalized Jacobian (Theorem 1.13 \cite{DM}),
hence there cannot possibly exist a universal Picard scheme over the whole of $\mgb$
(for the same reason why there exists no universal curve). 

The stack
theoretic approach is  thus necessary to answer the above question in general; 
the corresponding result  is the following:  there exists a smooth
Deligne-Mumford stack $\pdst$, with a natural representable morphism to 
the stack  $\mgbst$, such that 
for every
family $f:\X
\la
B$ of  stable curves 
with  $\X$  regular, the N\'eron model of $\Pic^d_K$ is the fiber product
$\pdst\times_{\mgbst}B$.

The stack $\pdst$ has a geometric description, as it corresponds to the 
 ``balanced Picard functor",
which is a separated partial completion of the 
 degree-$d$ component of the classical
Picard functor on smooth curves (cf. \ref{balfun} and \ref{modstack}). 
Similarly the  scheme $\pd$ is  the fine moduli space for such a functor
restricted to automorphism-free curves (\ref{fine}).

The requirement $(d-g+1,2g-2)=1$
is  well known (by \cite{MR})
 to   be necessary and sufficient for the existence of a Poincar\'e line
bundle for the universal Picard variety 
$\Pic^d_g\la M_g^0$
(associated to the universal family of smooth curves);
we extend such a result 
as follows.
Our scheme  $\pd$ will be constructed as a dense open subset of
the compactification $\pdg$ of $\Pic^d_g$
obtained in 
\cite{caporaso};
we prove that the above Poincar\'e line bundle extends 
over $\pdg$.
More precisely, we prove that such a numerical condition 
characterizes when the balanced Picard functor is representable (and separated), and 
  when the
corresponding groupoid is a Deligne-Mumford stack, 
representable over $\mgbst$ (cf. chapter~\ref{represent}).
Thus  the hypothesis that $(d-g+1,2g-2)=1$ plays a crucial role in various places
of our argument;
 we are therefore  led to conjecture that 
without it
the parametrization result 
(\ref{picner}) would fail.

A consequence of the construction is 
 a modular completion  $\pdgst$ of $\pdst$ by a smooth
Deligne-Mumford stack  representable over $\mgbst$, which enables us to
    obtain a 
geometrically meaningful compactification of
the N\'eron model for every family $f$ as above.

We prove  that our compactification of the N\'eron model is endowed with a 
canonical stratification described in terms of the
N\'eron models of the connected partial normalizations of the closed fiber of $f$
(Theorem~\ref{strata}).
Moreover, in \ref{quot}, we exhibit  it as a ``quotient" of the N\'eron
model for 
a 
ramified degree $2$-base change of $f$.

Notice that as $d$ varies, the  closed fibers of $\pd\la \mgb$  do not, 
hence the question naturally arises
as to  how many isomorphism classes of these spaces there are;
the exact number of them is computed in \ref{euler}.
\end{nota}

\begin{nota}{\it Context.}
\label{}
The language and the techniques used in this paper are mostly those of
\cite{BLR} for the theory of N\'eron models,
and of \cite{GIT} for Geometric Invariant Theory.

As we said, we use the compactification $\pdg \la \mgb$ of the universal Picard variety;
however
such a space existed   only   as a scheme,
 not as a stack.
To answer our initial question about the parametrization of N\'eron models,
we  need to
 ``stackify" such a  construction  and to build the 
standard universal elements for it (the universal curve and the
Poincar\'e bundle). 
This occupies most of section \ref{represent}.    
We are 
in a   lucky position to
  apply the  theory of stacks as was  
 developed in recent years, in fact 
$\pdg$ and $\pd$ are geometric GIT-quotients hence
our stacks are  ``quotient stacks",  which have been
carefully studied by many authors.
In particular, we use  \cite{AV}, \cite{ACV}, \cite{edidin}, \cite{LM}
and \cite{vistoli}, together with the seminal paper \cite{DM}.

Why should $\pdst $ be a good candidate to glue N\'eron models
together over $\mgbst$? The initial observation, already at the scheme level, is that
if the condition $(d-g+1,2g-2)=1$ holds  every closed fiber of
$\pdg$ over  $\mgb$
contains the fiber of the corresponding N\'eron model as a dense open subset. 

\

N\'eron models provide
 the solution for a fundamental mapping problem (see the ``N\'eron  mapping
property" in \ref{notner}) and are uniquely determined by this.
Their existence for abelian varieties was established by A. N\'eron
in \cite{neron}; the theory was  developed by M. Raynaud  (in \cite{raynaud})
who, in particular, unraveled the connection with  the Picard
functor in a way that will be heavily used in this paper.
N\'eron models have  been widely applied in arithmetic and algebraic
geometry; 
a remarkable example is
 the   proof (valid in all characteristics) of the 
stable reduction theorem for curves given in \cite{DM}.
Nevertheless they  rarely appear in the present-day moduli theory of
curves, where their potential impact looks  promising (see Section~\ref{abel}).

N\'eron models are well known not to have good functorial properties: 
their formation does not commute
with base change, unless it is an \'etale one. 
However there are
  advantages in having a geometric description for them  (and
 for their completion), such as the possibility to  interpret mappings  in a geometric way
(note that their universal property gives us the  existence of many such mappings,
some arising form remarkable geometric settings).
This may be fruitfully used to study problems concerning limits of line bundles and
linear series, as we  briefly illustrate in \ref{abel}.

We  mention
one further motivating issue;
that is
 the problem of comparing  various existing 
completions  of  the
Picard functor and of some of its distinguished subfunctors
(such as the spin-functors or the functor of torsion points
in the Jacobian).
It is fair to say that our understanding of the situation is insufficient,
a clear picture of how the various
 compactifications mentioned above relate
to each other is missing.
An overview of various completions of the generalized Jacobian with some comparison
results is in \cite{alexeev} (more details   in     \ref{many}); 
the interaction between
compactified   spin schemes and Picard schemes is studied
in
\cite{fontanari} and  \cite{CCC}; various basic questions remain open.
Understanding the relation with
N\'eron models   can be used for such problems,
thanks to the N\'eron   mapping property  (see \ref{compare}).
\end{nota}
\begin{nota}{\it 
Summary.}
\label{}
The paper is organized as follows: section \ref{twister} recalls some basic facts about our
N\'eron models,
sections \ref{picard} and 
\ref{represent} are about the 
``balanced Picard functor" and the corresponding stack;
in \ref{comparison}  
the  connection with N\'eron models
is established, together with some comments and examples. The last two sections are devoted
to the completion of the N\'eron model, which  is described in \ref{completion}
with focus on the stratification,
and in \ref{quotient} as a quotient of a 
N\'eron model of a certain base
change.  In the appendix some comments about applications, together with  some
useful combinatorial facts,  are collected.
\end{nota}

\noindent{\it Aknowledgments.}
I  wish to express my gratitude to Dan Abramovich and Angelo Vistoli
for  their  kind  explanations of crucial help for
section~\ref{represent}, and  to Cinzia Casagrande and Eduardo Esteves for various useful comments.

\section{Notation and terminology}
\label{not}

\begin{nota}
\label{notfun}
All  schemes are assumed  locally of finite type
over an an algebraically closed field
$k$, unless otherwise specified.
$R$ denotes a discrete valuation ring (a DVR) with algebraically closed residue field
$k$ and quotient field  $K$. 
For any scheme $T$ over  $\Spec R$ we  denote $T_K$ the generic fiber and $T_k$ the closed fiber.
 
If $\phi:W\la B$ is a morphism and $T\la B$ is a $B$-scheme
we shall denote $W_T:=W\times _BT$ and $\phi_T:W _T \la T$ the 
projection.
\end{nota}
\begin{nota}
\label{notsub} 
$X$ will be a nodal  connected curve projective over $k$;   
$C_1,\ldots ,C_{\gamma}$ its irreducible components. 

For any complete subcurve $Z\subset X$,  $g_Z$ is its arithmetic genus, 
$Z':=\overline{X\smallsetminus Z}$ its  complementary curve and 
$k_Z:=\#(Z\cap Z')$. Then
$$
w_Z:=\deg_Z\omega_X=2g_Z-2+k_Z.
$$
For a line bundle $L\in \Pic X$ its {\it multidegree}  is 
$ \mdeg L :=(\deg _{C_1}L,\ldots,\deg _{C_{\gamma}}L)$.

We denote
 $\md = (\vectd)$ elements of $\Z^{\gamma}$ 
(or in $\Q^{\gamma}$) and $|\md |:=\sum_{1}^{\gamma} d_i$.
We say that $\md$ is positive (similarly, non-negative, divisible by some integer, etc.) if all $d_i$ are.
If $\md \in \Z ^{\gamma}$
(or in  $\Q^{\gamma}$)  we denote the  ``restriction of $\md$ to the subcurve $Z$" 
of $X$ by
$
d_Z=\sum _{C_i\subset Z}d_i
$.

We set
$
\picX{\md}:=\{L\in \Pic X: \mdeg L = \md\}
$.
In particular, the {\it generalized Jacobian} of $X$ is 
$ \picX{\underline 0}:=\{L\in \Pic X: \mdeg L = (0,\ldots, 0)\}.$
There are (non-canonical) isomorphisms $\picX{\underline{0}} \cong \picX{\md}$
for every $\md \in \Z^{\gamma}$.
Finally we set
$
\picX{d}:=\{L\in \Pic X:\deg L = d\}= \coprod_{|\md|=d}\picX{\md}.
$

\end{nota}

\begin{nota}
\label{notpic}

$ f:\X \la B$ will denote
a family of  nodal curves;
that is, $f$ is a proper flat morphism 
of schemes  over $k$,
such that every closed fiber of $f$ is a
connected  nodal curve.

$
 \Picfun _f  
$ denotes
the Picard functor of such a family (often denoted 
$ \Picfun _{\X/B}$ in the literature, see \cite{BLR} chapter 8
 for the general theory).
$  \Picfun _f^d$ is
the 
subfunctor of line bundles of (relative) degree $d$.

We shall often consider $B=\Spec R$; in that case the closed fiber of $f$ will be
denoted by $X$;  let us assume
that for the rest of the section.
$\Picfun _f$ 
(and similarly $\pfun{d}$)
is represented by a scheme $\Pic_f$
(due to D. Mumford, see \cite {BLR} Theorem 2 in 8.2 and \cite{mumford})
which
may very well fail to be separated:
if all geometric fibers of $f$ are irreducible, then $\Pic_f$  is separated
(due to A. Grothendieck \cite{SGA}, see also \cite {BLR} Theorem 1 in 8.2) and conversely
(see
\ref{nonsep}).

The identity component of the Picard functor 
is well known to be represented by a separated scheme over $B$
(the  {\it generalized Jacobian}, see \cite{raynaud} 8.2.1), which we
shall denote 
$\jacf$
(denoted by $P^0$ in \cite{raynaud} and   by
$\Pic ^0_{\X/B}$in \cite{BLR}).
 
For any $\md \in \Z ^{\gamma}$ 
  consider 
$\picf{\md}\subset \Pic _f^d$,   parametrizing 
line bundles of degree $d$ whose restriction to the closed fiber has multidegree $\md$.
Just like
$\jacf$, these are fine moduli schemes;
 $\picf{\md}$ is a natural
$\jacf$-torsor .

The generic fiber of
$\jacf$ and of $\picf{0}$ coincide and will be denoted by $\jac$;
similarly $\jacd$ denotes the generic fiber of $\picf{d}$
(and of $\picf{\md}$).

\end{nota}
\begin{nota}
\label{notstab}

A {\it stable} curve is (as usual) a nodal connected curve of genus $g\geq 2$ having 
ample dualizing sheaf. The moduli scheme (respectively stack) for stable curves of genus $g$ is 
denoted by
$\mgb$ (resp. $\mgbst$). If $g\geq 3$ the locus $\mgb^0\subset \mgb$ 
of curves with trivial automorphism group is  nonempty, open and nonsingular.

A {\it semistable} curve is a nodal connected curve of genus $g\geq 2$ whose dualizing sheaf
has non-negative multidegree.
A {\it quasistable} curve $Y$ is a semistable curve 
such that any two of its  {\it exceptional components} do not meet
(an exceptional component of $Y$
is a smooth rational component $E\cong \pr{1}$  such that $\#(E\cap
\overline{Y\smallsetminus E})=2$).

If $Y$ is a semistable curve, its {\it stable model} is the stable curve obtained by contracting
all of the exceptional copmonents of $Y$. For a given stable curve  $X$ there exist  finitely many
quasistable curves having
$X$ as stable model; we shall call such curves the
{\it quasistable curves of $X$}.
\end{nota}
\begin{nota}
\label{notner}
Let $B$ be a connected Dedekind scheme  with function field $K$.
 If $A_K$ is an abelian variety over $K$,
or a torsor under a smooth group scheme, we denote by ${\operatorname N}
(A_K)$ the N\'eron model of $A_K$, which is a smooth  model of $A_K$
over
$B$ uniquely determined by the following universal property
(the {\it N\'eron mapping property}, cf. \cite{BLR} definition 1):
  every $K$-morphism
$u_K:Z_K\la A_K$ defined on  the generic fiber of some scheme 
$Z$ smooth over $B$ admits a unique extension  
to a $B$-morphism $u:Z\la N(A_K)$.

Recall that $N(A_K)$  may  fail to be proper over $B$,
whereas it is obviously separated.
Although $ N(A_K)$ is endowed with a canonical torsor structure, induced by the one of $A_K$, we
shall always consider  it merely as a scheme.
\end{nota}

\section{The N\'eron model 
for the  degree-$d$ Picard scheme}
\label{twister}

We begin by introducing  N\'eron models of Picard varieties of curves,
 following  Raynaud's approach (\cite{raynaud}).
Most of the material in 
this section   is in chapter 9 of \cite{BLR}
in a far more general form
( also in sections 2 and 3 of \cite{edix}, which is closer to our situation); we
revisit it with a slightly different terminology, suitable to our goals.

\begin{nota}
\label{nonsep} 
Let $f:\X \la B=\Spec R$ be a family  of  curves and denote
 $\gen$ the generic fiber, assumed to be smooth.
To  construct the N\'eron model of the Picard variety
$\Pic^d_K:=\Pic^d\gen$ 
of $\gen$,  it is  natural to look at the Picard scheme $\Pic_f^d\la B$
of the given family, which is smooth and has generic fiber equal to $\Pic^d_K$.
The problem is that  $\Pic_f^d$ will fail to be separated over $B$ as soon as the closed fiber
$X$ of $f$
 is  reducible.

From now on we  assume that $X$ is a reduced curve
 having at most nodes as singularities.
Its decomposition  into  irreducible components is denoted
$
X=\cup _1^{\gamma}C_i.
$
\end{nota}
One begins by isolating line bundles on $X$ that are specializations of the trivial bundle
(so called {\it twisters}).
\begin{defi}
\label{twidef} 
\begin{enumerate}[(i)]
\item
Let $f:\X\la \Spec R$ be a family of nodal curves.
A line bundle $T\in \Pic X$ is called an {\it $f$-twister}
(or simply a {\it twister}) if there exist integers $\vectn$ such that 
$
T\cong \O_\X(\sum_1^{\gamma} n_iC_i)\otimes \O_X
$
\item
The set of all $f$-twisters is a discrete subgroup of $\Pic^0X$, 
denoted $\twX$.
\item
Let $L, L'\in \Pic X$. We say that $L$ and $L'$ are
{\it $f$-twist equivalent} (or just {\it twist equivalent})
if for some $T\in \twX$  we have
$
L^{-1}\otimes L'\cong T
$
\end{enumerate}
\end{defi}
The point is: every separated completion of $\Pic{\gen}$
over
$B$ must identify twist equivalent line bundles.

\begin{remark}
\label{prifib}
Notice that the integers $\vectn$ are not uniquely determined, as 
$X$ is a principal divisor (the base being $\Spec R$) and
we have for every integer $n$,\  
$
\O_\X(nX)\otimes \O_X\cong \O_X
$
\end{remark}

We need the following well known 
(6.1.11 in \cite{raynaud} and
 \cite{BLR}, lemma 10. p. 272)
list of facts (recall that $k_Z:=\#Z\cap\overline{X\smallsetminus Z})$):

\begin{lemma}
\label{twilm}
Let $f:\X \la \Spec R$ be a family of nodal curves with $\X$  regular.
\begin{enumerate}[(i)]
\item
\label{twideg}
$
\mdeg (\O_\X(\sum_1^{\gamma} n_iC_i)\otimes \O_X)
=\underline{0}$
if and only if $n_i=n_j$ for all $i,j=1\ldots \gamma$.
\item
\label{twinb}
Let $T$ be a non-zero twister. There exists a subcurve $Z\subset X$
such that  $$\deg _ZT\geq k_Z$$
\item
\label{twiseq}
There is a natural exact sequence
$$
0\la \Z \la \Z^{\gamma} \la \twX \la 0.
$$
\end{enumerate}
\end{lemma}
\begin{proof}
For (\ref{twideg}), 
set
$T= \O_\X(\sum_1^{\gamma} n_iC_i)\otimes \O_X$.
One direction follows immediately from 
(\ref{prifib}).
Conversely
 assume that $\mdeg T=0$.
Define for $n\in \Z$ the subcurve $D_n$ of $X$ by
$$
D_n = \cup _{n_i=n}C_i
$$
(If $n=0$
the curve $D_0$  is the union of all components  having coefficient $n_i$ equal to zero.)
Now $X$ is partitioned as 
$
X=\cup _{n\in \Z}D_n
$
and every irreducible component of $X$ belongs to exactly one $D_n$. 
By construction
$$
T=\O_\X(\sum_{n\in \Z} nD_n)\otimes \O_X
$$
and our goal is to prove that there is only one nonempty $D_n$
appearing  above.
Let   $m$ be the minimum integer such that $D_m$ is not empty, thus 
$D_n =\emptyset$ for all $n<m$.
We have
\begin{eqnarray}
\label{twieq}
\lefteqn{deg _{D_m}T=-mk_{D_m}+\sum_{n>m} n(D_n\cdot D_m)\geq }\nonumber\\
&& -mk_{D_m}+(m+1)\sum_{n>m} (D_n\cdot D_m)\geq \sum_{n>m} (D_n\cdot D_m)=k_{D_m}\geq 0
\end{eqnarray}
where in the last inequality we have equality if and only if all $D_n$ are empty for $n>m$
(so that $X=D_m$). On the other hand the hypothesis was 
$\mdeg T=0$ and hence equality must hold above, so we are done.
This also proves (\ref{twinb}) by taking $D_m=Z$.

Now we prove (\ref{twiseq}).
The sequence is defined as follows
$$\begin{array}{lcccccccr}
0&\la& \Z  & \stackrel{\sigma}{\la} &\Z^{\gamma} &\stackrel{\tau}{\la} &\twX &\la &0 \\
& & 1 &\mapsto &(1,\ldots,1) &&&&\\
&&&&(\vectn )& \mapsto &\O_\X(\sum_1^{\gamma} n_iC_i)\otimes \O_X

\end{array}
$$
The map $\tau$ defines a Cartier divisor because $\X$ is regular.
The injectivity of $\sigma$ and the surjectivity of $\tau$ are obvious.
The fact that $\im \sigma \subset \ker \tau$ was oserved before
(in (\ref{prifib})). Finally, suppose that $(\vectn )$ is such that
the associated $f$-twister $T:=\O_\X(\sum_1^{\gamma} n_iC_i)\otimes \O_X$ is zero.
Then $T$ must  have  multidegree equal to zero, therefore, by the first part, we
obtain that $(\vectn )=(m,\ldots,m)$ for some fixed $m$ and hence $(\vectn )\in \im
\sigma$.
\end{proof}

\begin{nota}
\label{discrete}
Twisters on a  curve $X$ depend on two types of data:
(1) discrete data, i.e. the choice of the coefficients $\vectn$,
(2) continuous data, namely the choice of  $f:\X\la B=\Spec R$. 
More precisely,
while  twisters  may depend on $f:\X\la B$,
their multidegree 
only depends of the type of singularities of $\X$
(see \ref{type}). Let us assume that   $\X$ is
regular. For every component $C_i$ of $X$ denote, if $j \neq i$,\ 
$
k_{i,j}:= \#(C_i\cap C_j)
$
and 
$
k_{i,i}=-\#(C_i\cap \overline{C\setminus C_i})
$
so that the matrix
$
\int := (k_{i,j})
$
is an integer valued symmetric matrix which can be viewed as an {\it intersection matrix} for $X$.
It is clear that for every pair $i,j$ and for every (regular) $\X$,
$
\deg _{C_j}\O_\X(C_i)=k_{i,j}.
$
We have that 
$\sum _{j=1}^{\gamma} k_{i,j}=0$ for every fixed $i$.
Now, for every $i=1,\ldots ,\gamma$ set 
$
\mc :=(k_{1,i},\ldots ,k_{\gamma ,i})\in \Z ^{\gamma}
$\  
and
$$
{\bf Z}:=\{ \md \in \Z ^{\gamma}: \  |\md | = 0\}
$$
so that $\mc\in {\bf Z}$ and we can consider the sublattice $\LX $ of ${\bf Z}$ spanned by them 
$$
\LX :=< \underline{c}_1, \ldots ,\underline{c}_{\gamma} >.
$$
Thus,  $\LX$ is the set of multidegrees of all twisters 
and has rank $\gamma -1$
(by \ref{twilm}  (\ref{twiseq})).
\end{nota}
\begin{defi}
\label{dcg}
The {\it degree class group} of $X$ is the (finite) group
$\dcg := {\bf Z}/\LX$. 
Let $\md$ and $\md '$ be in $\Z^{\gamma}$; we say that they are equivalent,
denoting $\md \equiv \md '$,  iff
their difference is the multidegree of a twister, that is if
$\md - \md '\in \LX$
\end{defi}

\begin{nota}
\label{dcgref}
The degree class group 
is a natural invariant to consider in this setting,  it was first (to our knowledge)
defined and studied by Raynaud 
(in  8.1.2 of \cite{raynaud},  denoted  $\ker \beta/\im \alpha$). 
We here adopt  the terminology and notation used in
 \cite{caporaso} section 4.1, which is convenient for our goals.

$\dcg$   is  the 
component-group of the N\'eron model of the Jacobian of a family of nodal curves
$\X
\la \Spec R$ with $\X$ nonsingular (see thm.1 in 9.6 of \cite{BLR} and also \ref{dc}).
The group of components of a N\'eron model,
in more general situations than the one studied in this paper, 
has been the object of much research. In particular, bounds for its
cardinality have been obtained by D. Lorenzini in \cite{lorgroup}; see also
\cite{lorgraph}, \cite{lorner} and\cite{BL}   for further study 
and applications.
It is quite clear that $\dcg$ is  a purely combinatorial invariant of 
$X$, a description of it in terms of the dual graph (due to Oda and Seshadri \cite{OS})
is recalled in \ref{degcomp}. 
\end{nota}

\begin{nota}
\label{trans}
The  group $\dcg$ parametrizes classes of multidegrees summing to zero.
More generally, let us denote $\dcg^d$ the set of  classes of multidegrees 
summing to $d$:
$$
\dcg^d:= \{ \md \in \Z ^{\gamma}: |\md | =d\}/_\equiv
$$
(where ``$\equiv$" is defined in \ref{dcg}).
We shall denote the elements in $\dcg^d$ by lowercase greek letters $\delta$ and write
$\md \in \delta$ meaning that the class $[\md]$ of $\md $ is $\delta$. Of course,
there are bijections $\dcg^d\leftrightarrow \dcg$.
\end{nota}
\begin{nota}
\label{identify}
Let  $f:\X \la \Spec R=B$ with $\X$ regular and, as usual, assume that the closed fiber has $\gamma$ irreducible components. Let
$\md$ and
$\md '$ be  equivalent 
multidegrees, then there is a canonical isomorphism (depending only on $f$)
$$
\iota_f(\md,\md'):\Pic_f^{\md}\la \Pic_f^{\md'}
$$
which restricts to the identity on the generic fiber. To prove that, recall that by
\ref{twilm} part (\ref{twideg}) there exists a unique $T\in \twX$ such that
$
\mdeg T =\md '-\md$
and that there is a unique line bundle ${\mathcal T}\in \Pic \X$ such that
${\mathcal T}$ is trivial on the generic fiber and
${\mathcal T}\otimes \O _X \cong T$; in fact ${\mathcal T}$ must be of the form $\O_{\X}(\sum n_iC_i)$ and the $n_i$ are determined  
up to adding a multiple of the closed fiber
(see \ref{prifib}), which does not change the equivalence
class of ${\mathcal T}$, as $X$ is a principal divisor on ${\mathcal X}$ ($\Pic
B=0$). The isomorphism $\iota=\iota_f(\md,\md')$ is thus given by tensor product
 with ${\mathcal T}$, so that
if $L\in \Pic _k^{\md}$ we have $\iota (L)=L\otimes T$, whereas if $L\in \Pic^d_K$
then
$\iota(L)=L$.

We shall therefore  identify $\Pic _f^{\md}$ with $\Pic _f^{\md'}$ for all pairs of equivalent
$\md$,
$\md '$. Thus for every $\delta \in \dcg^d$ we define for every $\md \in \delta$
\begin{equation}
\label{picdel}
\picf{\delta}:=\picf{\md}\  \  
\end{equation}

The schemes  $\picf{\delta}$ for a fixed total degree $d$ all 
have the same generic fiber, $\Pic_K^d$; we can then glue them together
identifying their
 generic fibers. We shall denote the so obtained scheme over $B$
$$
 \frac{\coprod _{\delta \in \dcg^d}\picf{\delta}}{\sim}
$$
(where $\sim$ denotes the gluing along the generic fiber) so that its
  generic fiber  is
$\Pic^d_K$.
We have 
\end{nota}

\begin{lemma}
\label{sep}
Let $f:\X \la \Spec R$ be a family of nodal curves
with $\X$
regular. Then we have a canonical $B$-isomorphism
$$
\NN(\Pic^d_K)\cong \frac{\coprod _{\delta \in \dcg^d}\picf{\delta}}{\sim}.
$$
\end{lemma}
\begin{proof}
We may  replace $B$ by its strict henselization, in fact all the
objects involved in the statement are compatible with \'etale base changes
(of course $\X$ remains regular under any such  base change, and $\dcg$ 
does not change).
Recall also that N\'eron models descend from the strict henselization of $B$ to $B$ itself
(\cite{BLR} 6.5/3).

Assume first that $d=0$.
The N\'eron model of $\Pic ^0_K$ is proved in \cite{BLR} (Theorem 4 in 9.5)
to be equal to  the quotient
$\Pic^0_f/E$ where $E$ is the schematic closure of the unit section
$\Spec K \la \Pic ^0_K$
(so that $E$ is  a scheme over $B$, see \cite{BLR} p. 265).

We can explicitly describe the closed fiber of $E$:
$E_k=\twX$. In fact if $L$ belongs to the closed fiber of $E$, then $L$ 
is a line bundle on $X$ 
which is a specialization of the trivial line bundle on $\gen$;
thus
there exists a line bundle $\L$ on the total space $\X$ which is trivial on the generic fiber of $f$ and whose restriction to $X$ is $L$. Therefore $\L$ is of the form $\L =
\O_\X(D)$ with $D$ supported on $X$,
hence
$L\in \twX$. The converse, i.e. the fact that $\twX$ is in $E_k$, is obvious.
Now we have 
$$
\Pic ^0_f=\frac{\coprod_{|\md|=0}\Pic^{\md}_f}{\sim}
$$
where $\sim$ denotes (just as above) the gluing of the schemes $\Pic^{\md}_f$ along
their
 generic fiber (which is the same for all of them: $\Pic ^0_K$).

We obtain that the quotient by $E$  identifies $\Pic^{\md}_f $ with $\Pic^{\md'}_f$ for all
pairs
 of equivalent $\md$ and $\md'$, and this identification  is the same induced
by
 $\iota_f(\md, \md')$ which was used to define $\Pic_f^{\delta}$ in \ref{identify}
formula (\ref{picdel}). Hence we have canonical isomorphisms
$$
\Pic^0_f/E\cong\frac{\coprod_{|\md|=0}\Pic^{\md}_f}{\sim}\cong
\frac{\coprod _{\delta \in \dcg}\picf{\delta}}{\sim}.
$$

 For general $d$, we have that $\Pic ^d_K$ is a trivial
 $\Pic ^0_K$-torsor (in the sense of \cite{BLR} 6.4) 
and we can reason
as we just did to obtain
$$
\NN(\Pic^d_K) = \Pic^d_f/E^d=(\coprod_{|\md|=d}\Pic _f^{\md})/E^d\cong
(\coprod_ {\delta \in \dcg^d}\Pic _f^{\delta})/\sim 
$$
where  $E^d$ denotes the analog of $E$, that is the schematic closure
of a fixed section $\Spec K \la \Pic ^d_K$
(which exists because, $R$ being henselian,  $f$ has a section).
\end{proof}

\begin{remark}
\label{dc}
The lemma clarifies \ref{dcgref}:
the degree class group $\dcg$ is the  group of connected components
of 
the closed fiber of $\NN(\Pic^0_K)$.
In fact (recalling \ref{trans}) for
the closed fiber  we have  
$$
(\NN(\Pic^d_K))_k\cong\picX{d}/\twX \cong \coprod _{\delta \in \dcg^d}\picX{\delta}.
$$

\end{remark}

\section{The balanced Picard functor}
\label{picard}

  As stressed  in \ref{dc}, the scheme structure of the closed fiber of the N\'eron
model  does not depend on 
the family $f$ (the hypothesis that $\X$ is a nonsingular surface is crucial, see
\ref{type}). We shall now ask
whether, for a fixed
$d$, our N\'eron models ``glue together"  over $\mgb$.
From  the previous section, a good starting point would be to 
to find a ``natural" way of choosing representatives
for  multidegree classes.

\begin{example}
\label{motex}
Let $d=0$ and consider the identity in $\dcg$; then $(0,\ldots , 0)$ is a natural representative for
that.  It is then reasonable to choose  representatives for the other classes 
so that their entries have the smallest possible absolute value.

For example, let $X=C_1\cup C_2$ with $C_1\cap C_2=k$ and $k$ odd.
Then $\Delta _X \cong \Z/k\Z$ and our choice is:
$$
(0,0),(\pm 1,\mp 1),\ldots ,(\pm \frac{k-1}{2},\mp \frac{k-1}{2})
$$
Another natural case  is $d=2g-2$; here the class $[\mdeg \omega _X]$,
represented of course by $\mdeg \omega _X$,
plays the role of the identity. Therefore, as before, the other representatives should be
chosen as close to $\mdeg \omega _X$ as possible.
For $X$ as  above
 the  representatives would be (recalling that $w_{C_i}:=\deg_{C_i}\omega_X$)
$$
(w_{C_1},w_{C_2}),(w_{C_1}\pm 1,w_{C_2}\mp 1),\ldots ,(w_{C_1}\pm
\frac{k-1}{2},w_{C_2}\mp
\frac{k-1}{2}).
$$
\end{example}
In what follows we use the notation of \ref{notsub}.
\begin{defi}
\label{dom} Let $X$ be a nodal curve of any genus.
\begin{enumerate}[(i)]
\item
\label{}

The {\it basic domain} of $X$ is the bounded subset $\domX\subset \Z^{\gamma}$
made of all $\md \in \Z^{\gamma}$ such that $|\md|=0$ and such that for every subcurve $Z\subset X$
we have
$$
-\frac{k_Z}{2}\leq d_Z\leq \frac{k_Z}{2}.
$$

\item
\label{}

For any $\mb \in \Q^{\gamma}$ such that $b:=|\mb|\in \Z$
denote $\domX(\mb)$ the subset of $ \Z^{\gamma}$
made of all $\md \in \Z^{\gamma}$ such that $|\md|=b$ and such that for every subcurve $Z\subset X$
we have
$$
b_Z-\frac{k_Z}{2}\leq d_Z\leq \frac{k_Z}{2} +b_Z
$$
\end{enumerate}
\end{defi}

\begin{remark}
\label{BIrmc}
Note that $\domX$ 
(and similarly $\domX(\mb)$)
is  the set of integral points contained in a polytope of
$\Q^{\gamma}$, whose boundary is defined by the inequalities in \ref{dom}.
We shall refer to $\domX(\mb)$ as a {\it translate} of $\domX$, although this is
is slightly abusive.
 
In the definition one could replace ``every subcurve $Z$ of $X$"
with ``every connected subcurve $Z$ of $X$" but not with ``every irreducible component of $X$".
In other words the basic polytope of $X$ is not in general, a product of $\gamma -1$
intervals (it is, of
course, if $X$ has only two components, in which case it is an interval).
\end{remark}

To connect with the previous discussion, we have
\begin{lemma}
\label{domlm}
Let $X$ be a nodal (connected) curve of any genus.
Fix any $\mb \in \Q^{\gamma}$ with $b:=|\mb|\in \Z$.
Then  every $\delta \in \dcg ^b$ has a representative contained in $\domX(\mb)$.
\end{lemma}
\begin{proof}
The proof
of proposition 4.1 in \cite{caporaso},  
apparently only a special case of this lemma (namely  $X$  quasistable (cf. \ref {notstab})
and $\mb = \mb_X^d$  as in (\ref{mb}) below),
carries out word for
word. 
\end{proof}
\begin{nota}
\label{domb}
We shall   choose a special  translate of $\domX^d$,
according to the topological characters of $X$.
Let $g\geq 2$,
set 
\begin{equation}
\label{mb}
\mb_X^d:=(w_{C_1}\frac{d}{2g-2},\ldots,w_{C_\gamma}\frac{d}{2g-2})\  \  \text{ and  }\  \  
\domX^d:=\domX(\mb_X^d)
\end{equation}
\end{nota}
then:
\begin{defi}
\label{BI}
Let $X$ be a semistable curve of genus $g\geq 3$ and $L\in\Pic ^dX$. Let $\md$ be the multidegree
 of $L$, We shall say that 
\begin{enumerate}[(i)]
\item
\label{semidef}
$\md$ is {\em { semibalanced}}
if  
for every subcurve $Z$ of $X$ the following 
(``Basic Inequality)" holds
\begin{equation}
\label{BIeq}
m_Z(d):=d\frac{w_Z}{2g-2}-\frac{k_Z}{2}\leq \deg_Z L\leq d\frac{w_Z}{2g-2}+\frac{k_Z}{2}=:M_Z(d)
\end{equation}
(equivalently, if $\md \in \domX^d$) and if for every exceptional component $E$ of $X$
\begin{equation}
\label{balexc}
0\leq \deg_E L\leq 1(=M_E(d)).
\end{equation}
\item
\label{semidef}
 $\md$ is {\em{balanced}} 
 if it is semibalanced and if for every exceptional component $E\subset X$
\begin{equation}
\label{balexc}
\deg_E L=1.
\end{equation}
\item
$\md$ is {\em{stably balanced}} if  it is
balanced and if for every subcurve $Z$ of $X$ such that
$
d_Z=m_Z(d)
$
we have that $\overline{X\smallsetminus Z}$ is a union of exceptional components.

\end{enumerate}
If $\X \la B$ is a family of semistable curves and $\L \in \Pic\X$ of reative degree
$d$, we say that
$\L$ is   (respectively {\em stably, semi}) {\em {balanced}} if for every $b\in
B$ the restriction of
$\L$ to
$X_b$ has (stably, semi) balanced multidegree.

\end{defi}
\begin{nota}
\label{dombal} 
In particular if $X$ is a stable curve the set $\domX^d$ (cf. \ref{domb})
equals the set of balanced multidegrees of total degree $d$.

The  inequality~(\ref{BIeq}) was discovered   by D. Gieseker 
in the course of the construction of
the moduli scheme
$\mgb$. Proposition 1.0.11 in \cite{gieseker} states that (\ref{BIeq}) 
is a necessary condition
for the GIT-semistability of the Hilbert point of a (certain type of) projective curve;
it was later proved in \cite{caporaso} that it is also sufficient. 
We mention that there exist other interesting incarnations of that inequality,
for example in \cite{OS} and \cite{simpson} (\cite{alexeev} connects them one to
the other). The
terminology  used in the above
definition was introduced in \cite{CCC}  (see  Theorem 5.16 there) to reflect the GIT-behaviour
of Hilbert points .
\end{nota}
 \begin{example}
\label{motexver}
The representatives in  \ref{motex} 
(for $d=0$ and $d=2g-2$) are all stably balanced and they are all the balanced
multidegrees for that $X$ and those $d$'s.
\end{example}

\begin{remark}
\label{semirk}
It is easy to check (combining (\ref{BIeq}) and (\ref{balexc}) of \ref{BI})
that  balanced line bundles live on
quasistable, rather than semistable curves, and hence on a ``bounded" class of curves. 
In analogy with  semistable curves, while 
semibalanced line bundles do
not admit a nice moduli space  (just like semistable curves) they do admit a
``balanced line bundle  model"
(by contracting all of the exceptional components where the degree is $0$, see
\ref{sembal}).
\end{remark}
\begin{remark}
\label{va} Assume that $d$
is very large with respect to $g$, 
then a balanced line bundle $L$ on a quasistable curve $X$ of genus
$g$ is necessarily very ample. In fact 
if $Z\subset X$, it suffices to show that the restriction of $L$ 
to $Z$ is very ample;
 if $Z$ is exceptional then $\deg_ZL=1$, otherwise we have
$\deg_ZL\geq m_Z(d) =d\frac{w_Z}{2g-2}-\frac{k_Z}{2}$
and, since $w_Z\geq 1$ and $k_Z\leq g+1$, the claim follows trivially.
\end{remark}

\begin{remark}
\label{BIrm}
Notation as in \ref{BI}.
\begin{enumerate}[(a)]
\item
\label{BIext}
Set $Z':=\overline{X\smallsetminus Z}$.
Then $d=m_Z(d)+M_{Z'}(d)$, in particular $d_Z=m_Z(d)$ if and only if $d_{Z'}=M_{Z'}(d)$.
\item
\label{BIst}
Let $X$ be  stable;  then $\md$ is stably balanced if and only if
strict inequality holds in (\ref{BIeq}) for every $Z$.
\item
\label{BIq}
Let $X$ be  quasistable. Then  a balanced
$\md$ is stably balanced if and only if the  subcurves  where strict inequality
in (\ref{BIeq}) fails are
all the $Z'$ unions  of exceptional components (in which case $d_{Z'}=M_{Z'}(d)$)
and (by (\ref{BIext})) their complementary curves $Z$
(in which case $d_Z=m_Z(d)$).
\end{enumerate}
\end{remark}

\begin{prop} 
\label{balrep}
Fix $d\in \Z$ and  $g\geq 2$.
\begin{enumerate}[(i)]
\item
\label{balrepe}
Let $X$ be a  quasistable curve of genus $g$ and $\delta \in \dcg^d$.
Then $\delta$ admits a semibalanced representative.
\item
\label{balun} A balanced multidegree is unique in its equivalence class if
and only if it is stably balanced.

\item
\label{balrepn} $(d-g+1,2g-2)=1$
if and only if for every quasistable curve $X$ of genus $g$ and  every $\delta \in
\dcg^d$,   $\delta$ has a  unique semibalanced representative.
\end{enumerate}
\end{prop}

\begin{proof} (\ref{balrepe})
By \ref{domlm} we know that every $\delta$ has a representative $\md$ in $\domX^d$; 
if $X$ is stable this
is enough. Assume that 
$X$ has an exceptional component $E$, notice that $m_E(d)=-1$ thus
we must prove that a representative for $\delta$ can
be chosen so that its restriction to $E$ is not $-1$. Assume first that $E$ is the unique exceptional
component.
 Observe that
for any subcurve
$Z\subset X$ and every decomposition 
$Z=A\cup B$ into
two subcurves having no component in common and meeting in $k_{A,B}$ points, we have
(omitting the dependence on $d$ to simplify the notation)
\begin{equation}
\label{maxeq}
M_Z=M_A+M_B-k_{A,B}.
\end{equation}
Now let $\md \in \domX^d$ and suppose that $d_E=-1=m_E$, denote $Z=E'$ the complementary curve and 
note that by  \ref{BIrm} (\ref{BIext}) we have that
$d_Z=M_Z$. Let $\underline{e}\in \LX$ be the multidegree  associated to $E$ (notation of
\ref{discrete}), we claim that $\md ':=\md - \underline{e}$ is semibalanced. 
We have that $d'_E=(\md - \underline{e})_E=-1+2=1$ so we are OK on $E$, now 
it suffices check every
connected subcurve $A\subset Z$ which meets $E$.  Suppose first that $E\not\subset A$, then
$d'_A=d_A - k_{A,E}$,
where $k_{A,E}=\#(A\cap E)>0$.
By contradiction assume that $\md'$ violates (\ref{BIeq}) on $A$, then,
as $d_A$ satisfies (\ref{BIeq}) and $d'_A< d_A$ we must have that 
$$d'_A=d_A-k_{A,E}<m_A=M_A-k_A$$
Now let $Z=A\cup B$ as above, so that
$k_{A,B}=k_A-k_{A,E}$, hence
$$d_A<M_A-k_{A,B}.$$
We conclude with the inequality
$$
M_Z=d_Z=d_A+d_B<M_A-k_{A,B}+M_B
$$
contradicting (\ref{maxeq}). 
Now let $E\subset A$; if $E\cap B = \emptyset$ then $d_A=d'_A$ and we are done. Otherwise $E$ meets $A$
in one point and one easily cheks that the basic inequality for $A$ is exactly the same as  for
$\overline{A\smallsetminus E}$, so we are done by the previous argument.

Since $X$ is quasistable, two of its exceptional components do not meet and
hence this argument can be iterated; this proves 
(\ref{balrepe}).

For (\ref{balun}), begin with a simple observation. 
 For every subcurve $Z$ of $X$,
the interval allowed by the basic inequality
 contains at most
$k_Z+1$ integers and the maximum
$k_Z+1$ is attained if and only if its extremes 
$m_Z(d)$ and $M_Z(d)$  are integers.   

Let now $\md$ be stably
balanced and $\underline{t}\in \LX$ (that is, $\underline{t}=\mdeg T$ for some twister $T$);
then, by
\ref{twilm} part (\ref{twinb}) there exists a subcurve $Z\subset X$  on which 
$(\md  + \underline{t})_Z\geq d _Z + k_Z$. This implies that  $\md  + \underline{t}$
violates the Basic Inequality, in fact
either  $d _Z$ lies in the interior of the allowed range and hence $d_Z + k_Z$ is out of the
allowed range;
or $d _Z$ is extremal, and   we use
\ref{BIrm}(\ref{BIq}).
 Therefore a stably balanced representative is
unique. Conversely, by what we said, two equivalent multidegrees that are both balanced
must be at the extremes of the allowed range for some curve $Z$, so neither
 can be stably balanced (by
\ref{BIrm}).

Now part (\ref{balrepn}).
As explained above,  it suffices to prove that
$(d-g+1,2g-2)=1$ if and only if for every $X$ quasistable of genus $g$ and every  subcurve
$Z\subset X$
such that neither $Z$ nor $Z'$ is  a union of
exceptional components,
$m_Z(d)$  is not integer.
Suppose that 
$(d-g+1,2g-2)=1$ holds, then $(d,g-1)=1$ (the converse holds only for odd $g$).
By contradiction, let
$X$ be a quasistable curve having a subcurve $Z$ as above for which $m_Z(d)$ is integer;
thus
\begin{equation}
\label{BIint}
\frac{dw_Z}{2g-2} = \frac{n}{2}\  \  \text{ where }\  \  n\in \Z :\  \  n\equiv k_Z   \mod (2)
\end{equation}
hence 
$g-1$ divides $w_Z$.
Then (by \ref{BIrm} (\ref{BIext}))
$M_{Z'}$ and $m_{Z'}$  are also  integer, therefore arguing as for $Z$,
$g-1$ divides $w_{Z'}$.
Now notice that
$
2(g-1)=w_Z+w_{Z'}
$
and  that $w_Z$ and $w_{Z'}$ are not zero (because $Z$ and $Z'$ are not union of exceptional
components). We conclude that
\begin{equation}
\label{BIg}
g-1=w_Z=w_{Z'} \  \  \text{ so that }\  \ g=2g_Z+k_Z-1.
\end{equation}
Thus  by the (\ref{BIint})
$$
\frac{dw_Z}{g-1} = d=n\  \  \text{hence }\  \  \  d\equiv k_Z   \mod (2).
$$
On the other hand  the second identity in (\ref{BIg}) shows that
$$(g-1)\equiv k_Z   \mod (2),   \  \  \text{ hence }\  \  \ d\equiv (g-1)   \mod (2).
$$
The latter implies that $2$ divides $(d-g+1,2g-2)$,   a contradiction.

Conversely: suppose that for some  $X$ and $Z\subset X$ we have
(see (\ref{BIint}))
$$
\frac{dw_Z}{g-1} = n\  \  \text{ with }\  \  n\in \Z :\  \  n\equiv k_Z   \mod (2)
$$
If $(d,g-1)\neq 1$ a fortiori $(d-g+1,2g-2)\neq 1$. 
Suppose then that $g-1$ divides $w_Z$;
we have just  proved that this implies $g-1=w_Z$, that $d=n$  and that 
$$
d\equiv (g-1) \mod (2)
$$
hence $2$ divides $(d-g+1,2g-2)$, and we are done.
\end{proof}

A weaker version of this result is proved in \cite{caporaso} sec. 4.2, where
the assumption that $d$ be very large is used. Despite the overlapping, we
gave here the full general
proof   to
stress the intrinsic nature of definition~\ref{BI} and contrast  the impression,
which may arise  from \cite{gieseker} and \cite{caporaso}, that it be a
technical condition deriving from Geometric Invariant Theory. 

A consequence of  \ref{balrep} and its proof is the following useful
\begin{cordef}
\label{dgen}
Let $d$ be an integer and $X$ a stable curve, we shall say that $X$ is $d$-{\em general}
(or  {\em general
for} $d$)
if  the following equivalent conditions hold.
\begin{enumerate}[(i)]
\item
\label{}
A multidegree on $X$ is balanced if and only if it is stably balanced.
\item
\label{}
The natural map  sending a balanced multidegree to its class
$$
\domX^d \la \dcg^d,\hskip.3in \md \mapsto [\md]
$$
 is a bijection.
\item
\label{}
For every  quasistable curve $Y$ of $X$, every element in $\Delta_Y^d$ has a unique
semibalanced representative.
\end{enumerate}
\end{cordef}

\begin{remark}
\label{g-1rm}
The assumption $(d-g+1,2g-2)=1$ in part
(\ref{balrepn}) of \ref{balrep} is a uniform condition 
ensuring that    every stable curve of
genus $g$ is $d$-general. The terminology is justified by the fact that the locus 
in $\mgb$ of $d$-general curves is open (see \ref{nofine}).

At the opposite extreme is the case $d=(g-1)$ 
(and, more generally, $d=n(g-1)$ with   $n$ odd), 
which is uniformly degenerate in the sense that for evey $X\in \mgb$ there exists
$\delta \in
\dcg^d$ having more than one balanced representative.
\end{remark}

We shall now define the moduli functor for balanced line bundles on stable curves.

\begin{defi}
\label{balfun}
Let $f:\X \la B$ be a family of stable curves and $d$  an integer. 
The 
{\it balanced Picard functor} $\pffun$
is  the contravariant functor from the category of $B$-schemes to the category of sets
which associates to a $B$-scheme $T$ the set of equivalence classes of balanced line bundles
$\L \in \Pic \X _T$ of relative degree $d$.
We say that $\L$ and $\L '$ are equivalent if there exists $M\in \Pic T$ such that
$
\L \cong \L ' \otimes f_T^*M.
$

A  $B$-morphism $\phi:T'\la T$ is mapped by $\pffun$ to the usual pull-back morphism
from $\pffun(T)$ to $\pffun(T')$.
\end{defi}

It is clear that $\pffun$ is a subfunctor of 
$\pfun{d}$. The point is that, in some  ``good" cases, $\pffun$ is 
representable by a separated scheme.

\begin{example}
\label{bag} Consider the   ``universal family of stable curves" of genus $g$
$$
f_g:\cgbar \la \mgb^0\subset \mgb .
$$ 
 (cf. \ref{notstab}). In this case we shall simplify the notation and set
$$
\pdfun:={\mathcal
 P}^d_{f_g}
$$
Observe that 
if $\pdfun$ is representable by a separated scheme $\pd$, then
for every family of automorphism-free stable curves $f:\X \la B$, the functor 
$\pffun$ is representable by the scheme
$\mu _f ^*\pd = B\times _{\mgb}\pd$
where $\mu_f:B\la \mgb$ is the moduli morphism of $f$.
\end{example}

\section{Balanced Picard schemes and stacks}
\label{represent}
The purpose of this section is to build the`` representable stack version" of the compactified
universal Picard variety constructed in \cite{caporaso} simply as a coarse moduli scheme.
\begin{nota}
\label{pdg}
From now  we fix integers  $d$ and $g\geq 3$  and we set $r:=d-g$.
We begin by recalling some    facts about the restriction of the balanced Picard functor
$\pdfun$ (cf. \ref{bag}) to
nonsingular curves (which is the ordinary Picard functor).

 The  degree-$d$ Picard functor for the universal family of nonsingular curves of
genus
$g$ is  denoted by $\Picfun^d_g$;  the so called ``universal degree-$d$ Picard
variety" over the moduli scheme of nonsingular  curves $M_g$ is
$\Pic _g^d\la M_g$ (we use here the notation ``$\Pic _g^d$" in place of 
``$P_{d,g}$" used in \cite{caporaso} and in \cite{HM}). The existence of the variety
$\Pic _g^d$ (a coarse moduli space in general, see below)
follows from  general results of A. Grothendieck 
(\cite{SGA} and \cite{mumford}, see \cite{GIT}
0.5 (d) for a summary).

Recall that,
for an arbitrary value of $d$, $\Picfun^d_g$
is only coarsely represented by $\Pic _g^d$, in fact a Poincar\'e
line bundle does not always exist. 
It is a well known
result due to N. Mestrano and S. Ramanan 
that $\Picfun^d_g$  is representable if and only if  $(d-g+1,2g-2)=1$ (in char$k=0$,
see Cor. 2.9 of
\cite{MR}). 
\end{nota}
\begin{nota}
\label{cap} We shall use 
the compactification 
$\pdg \la \mgb$
of  $\Pic _g^d\la M_g$  constructed in
\cite{caporaso}, from which   we need to recall and improve some results. Assume that
$d$ is very large (which is irrelevant, see below);
such a compactification is the GIT-quotient
 $\pdg=H_d/G$ 
 of the  action of the group
$G=PGL(r+1)$ on the
locus $H_d$ of
GIT-semistable points in the Hilbert scheme $\Hilb$
(for technical reasons concerning  linearizations, one actually carries out the GIT-construction
using the group $SL(r+1)$, rather than $PGL(r+1)$; 
since the two groups have the same orbits this will
not be not a problem).
\begin{enumerate}[{\bf (1)}]
\item
\label{capgen}
Denote by ${\mathcal Z_d}$ the restriction to $H_d$ of the universal family
over the Hilbert scheme
$$
\pr{r} \times H_d \supset{\mathcal Z_d}\la H_d,
$$
for $h\in H_d$ let $Z_h$ be the fiber of ${\mathcal Z_d}$ over $h$ and $L_h=\O_{Z_h}(1)$ the
embedding line bundle.
 $Z_h$
is a  nondegenerate quasistable  curve in $\pr{r}$ and $L_h$ is balanced 
in the sense of \ref{BI} (by \cite{gieseker});
conversely, every such a curve embedded  in $\pr{r}$ by a balanced line bundle 
appears as a fiber over $H_d$ (by \cite{caporaso}).
The point $h$ is GIT-stable  if and only if $L_h$ is stably balanced.
\item
\label{capstab}
For $h\in H_d$  denote $X_h\in \mgb$ the stable model of the quasistable curve $Z_h$. If $X_h$ is
$d$-general
 (see \ref{dgen}) the point $h$ is GIT-stable, which in turn implies that there
is a natural injection  (by \cite{caporaso} Section 8.2)
$$\Stab_G(h)\ha
\Aut (X_h).$$
Conversely, if $X\in \mgb$ is not $d$-general, there exists a
(strictly semistable)  $h\in H_d$ lying over $X$ having $\dim \Stab_G(h)>0$.
\item
\label{capreg}
 $H_d$ is regular and irreducible (by Lemma
2.2 and Lemma 6.2 in
\cite{caporaso}). 
\item
\label{capgeo}
The GIT-quotient $H_d/G$ is geometric (i.e. all  semistable points are 
stable)
if and only if $d$ is such that $(d-g+1,2g-2)=1$ (\cite{caporaso} Prop. 6.2).
\item
\label{caphigh}
For every pair of integers
$d$ and $d'$ such that
$
d\pm d' = n(2g-2),
$
 for  $n\in \Z$,
there are natural isomorphisms
$
\pdg \cong \overline{P}_{d',g}
$
(\cite{caporaso} Lemma 8.1). This allows us to define $\pdg$ for every $d\in \Z$,
  compatibly with the geometric description. That is, for $d\in
\Z$, pick $n$ such that $d'=d+n(2g-2)$ is large enough, the above isomorphism 
$
\pdg \cong \overline{P}_{d',g}
$ is constructed by tensoring with the $n$-th power of the relative dualizing sheaf. It is easy
to verify that a line bundle $L$ on a curve $X$ is balanced if and only if
$L\otimes \omega_{X}^{\otimes n}$ is balanced.
\end{enumerate}
\end{nota}

We begin with a scheme-theoretic result that will be generalized later on.
\begin{prop}
\label{fine} Let $g\geq 3$ and  $d$ be  such that $(d-g+1,2g-2)=1$.
\begin{enumerate}[(i)]
\item
The functor $\pdfun$ is representable by a separated scheme $\pd$.
\item
\label{}
 $\pd$ is integral,  regular and quasiprojective.
\item
\label{}
Let $[X]\in \mgb^0$ and denote $\pX$  the fiber of $\pd$ over it. Then $\pX$ is
regular of pure dimension
$g$.
In particular $\pd$ is smooth over $\mgb^0$.
\end{enumerate}
\end{prop}
\begin{proof}
Assume first that $d$ is very large ($d\geq 20(g-1)$ will suffice).
We use the notation and set up of \ref{cap} above. Denote by
$H_d^{st}$ the open subset of
$H_d$  parametrizing points corresponding to  stable
curves, that is
$$
H_d^{st}:=\{h\in H_d: Z_h \  \text{is a stable curve}\}.
$$
By \ref{cap} (\ref{capgen}) there is a natural surjective map $\mu :H_d^{st}\la \mgb$; set
$H:=\mu^{-1}(\mgb^0)$ so that
 $H$ parametrises points $h$ such that    $Z_h$ is a projective stable curve free from
automorphisms, 
$L_h$ is a degree-$d$ stably balanced
line bundle on $Z_h$ 
(by \ref{cap} (\ref{capgeo}))
and $\Stab_G(h)\cong \Aut(Z_h)=\{ 1\}$ (by \ref{cap} (\ref{capstab}))

We have that $H$ and $H_d^{st}$ are $G$-invariant integral
nonsingular  schemes (by \ref{cap} (\ref{capreg})). We shall denote 
$f_H:{\mathcal Z}\la H$ the restriction to $H$ of the universal family ${\mathcal Z}_d$
and define $\pd :=H/G$,
so that $H\la \pd$ is the geometric quotient 
of a free action of $G$.
Moreover, $G$ acts naturally (and freely) also on ${\mathcal Z}$ so that the quotent
${\mathcal C}_{\pd}:={\mathcal Z}/G$ gives a universal family on $\pd$.
Let us represent our parameter schemes and their  families
in a diagram
\begin{equation}
\label{finediag}
\begin{array}{lccccccr}
\pr{r}\stackrel{\pi}{\longleftarrow} \pr{r} \times H \supset&{\mathcal Z}&\stackrel{q}{\la} &
{\mathcal C}_{\pd}&\stackrel{p}{\la} &\cgbar &\\ 
&\  \  \downarrow{\tiny{f_H}}&&\downarrow&&\downarrow&\\
&H &\la & \pd=H/G &\stackrel{\phi}{\la}&\mgb^0&\subset\  \mgb
\end{array}
\end{equation}
Notice that all squares  are cartesian (i.e. fiber products)
so that all verticall arrows are universal families. 

Now let us consider the natural polarization
 $\L :=\O_{\mathcal Z}(1)=\pi ^* \O_{\pr{r}}(1)\otimes \O_{\mathcal Z}$.
As we said in \ref{cap},
$\L$ is stably balanced and, conversely,  
every pair $(X,L)$,  $X$  an automorphism free stable curve and $L\in \Pic ^dX$
a stably balanced line bundle, is represented by a $G$-orbit
in $H$.
More generally,
in Prop. 8.1 (2) of \cite{caporaso}  it is proved that $\pdg$ is a coarse moduli scheme for the
functor of stably balanced line bundles
on quasistable curves. 

In diagram (\ref{finediag}) we have exhibited  a universal family ${\mathcal C}_{\pd}\la \pd$, to
complete the statement we must show that
 there exists a universal 
or  {\it Poincar\'e line bundle} $\poin$ over ${\mathcal C}_{\pd}$
(determined, of course,  modulo pull-backs of line bundles on $\pd$).
This follows from lemma~\ref{poincare}, with $T=\pd$, $E=H$ and $\psi$ the inclusion,
so that $\X={\mathcal C}_{\pd}$.

We have so far proved that, if $d$ is  large,
the functor $\pdfun$ is represented by the scheme $\pd$ equipped with the universal pair
$({\mathcal C}_{\pd},\poin)$. 
The same result for all $d$ is obtained  easily
using \ref{cap} (\ref{caphigh}).

Now we prove (ii) and (iii).
We constructed $\pd$ as the quotient $H/G$
obtained by restricting the quotient $\pdg = H_d/G$,
that is, we have a diagram
\begin{equation}
\begin{array}{lccr}
H & \subset &H_d \\
\downarrow&&\downarrow&\\
\pd &\subset &\pdg .
\end{array}
\end{equation} 
Thus $\pd$ is quasiprojective because $H$ is open and $G$-invariant.
$\pd$ is integral and regular  because
$H$ is
irreducible and regular 
(\ref{cap} (\ref{capreg})) and $G$ acts freely on it. This concludes
 the second part of the statement.

The fact that 
 $\pX$ 
is smooth of pure dimension $g$ follows immediately 
from Cor. 5.1 in \cite{caporaso}, which implies that  
$\pX$ is a finite disjoint union of isomorphic copies of the generalized Jacobian of
$X$. 

Finally, $\pd$ is flat over  $\mgb^0$ 
(a consequence of the equidimensionality of the fibers) and, moreover, smooth
because the fibers are all regular.
\end{proof}
\begin{nota}
\label{notpoin}
Some notation before establishing the existence of Poincar\'e line bundles
and thus complete the proof of \ref{fine}.
If $\psi:E\la  H_d$ is any map we denote by $f_E:\ZZ_E=\ZZ_d\times _{H_d}E\la E$  and by
$\L_E\in \Pic \ZZ_E$ the pull back of the  polarization $\O_{\ZZ_d}(1)$
on $\ZZ_d$, so that $\L_E$ 
is a balanced line bundle of relative degree $d$.
 If, furthermore, $\pi:E\la T$ is a principal $G$-bundle and the above map $\psi$ is
$G$-equivariant, we can form the quotient
\begin{equation}
\begin{array}{lccr}
 \  \  \  \  \ZZ_E&\la &E \\
\  \  \  \  \  \downarrow &
 &\downarrow\\
\X=\ZZ_E/G&\stackrel{f}{\la}& E/G=T
\end{array}
\end{equation}
so that $f$ is a family of quasistable curves.

The proof of the next Lemma applies  a well known method  of M. Maruyama \cite{maruyama};
we shall make the simplifying assumption that $d$ be large, which will  later be removed.
\end{nota}
\begin{lemma}
\label{poincare} Notation as in \ref{notpoin}. Assume $d\gg 0$  and  $(d-g+1,2g-2)=1$.
Let $\pi:E\la T$ be a principal
$PGL(r+1)$-bundle  and  $\psi:E\la H_d$ be an equivariant map.
Then there exists a balanced line bundle $\poin \in \Pic \X$ of relative degree $d$ such that
for every $e\in E$ we have
$ (\L_E)_{|Z _e} \cong\poin_{|X _{\pi(e)}}$.
\end{lemma}
\begin{proof}
The statement holds locally over $T$, since $E\la T$ is a  $PGL(r+1)$-torsor.
Thus we can cover $T$ by open  subsets $T=\cup U_i$ such that, denoting 
the restriction of $f$ to
$\X_i:=f^{-1}(U_i)$ by
$$
f_i:\X_i\la U_i,
$$ there
is an $\poin _i\in \Pic
\X_i$ for which the thesis holds.
We  now prove that the $\poin _i$ can be glued together
to a line bundle over the whole of $\X$, modulo tensoring each of them by the pull-back of a line
bundle on
$U_i$.

By hypothesis there exist integers $a$ and $b$ such that
$$
a(d-g+1)+b(2g-2)=-1
$$
which we re-write as
\begin{equation}
\label{rew}
(a-b)(d-g+1)+b(d+2g-2-g+1)=-1
\end{equation}
Observe that, denoting by $\chi _{f_i}$ the relative Euler characteristic 
(with respect to the family $f_i$) we have that
$\chi _{f_i}(\poin_i)=d-g+1$ and $\chi _{f_i}(\poin_i\otimes \omega _{f_i})=d+2g-2-g+1$. 
Note also that
$\poin_i$ and $\poin_i\otimes \omega _{f_i}$ have no higher cohomology ($d$ is very large) and hence their direct images
via $f_i$ are locally free of rank equal to their relative Euler characteristic.
Define now for every $i$
$$
\N_i:=f_i^*\Bigl(\det (f_{i*}\poin_i)^{\otimes a-b}\otimes \det (f_{i*}\poin_i\otimes \omega
_{f_i})^{\otimes b}\Bigr).
$$

Now look at the restrictions of the $\poin _i$'s to the intersections $\X_i\cap \X_j$, we oviously have isomorphisms
$
\epsilon_{i,j}:(\poin _i)_{|\X_i\cap \X_j}\stackrel{\cong}{\la} (\poin _j)_{|\X_i\cap \X_j}
$
and hence for every triple of indeces $i,j,k$ an automorphism 
$$
\alpha_{ijk}:
(\poin _i)_{|\X_i\cap \X_j\cap\X_k}\stackrel{\cong}{\la} (\poin _i)_{|\X_i\cap \X_j\cap\X_k}
$$
where $\alpha_{ijk}=\epsilon_{k,i}\epsilon_{j,k}\epsilon_{i,j}$; thus
  $\alpha_{ijk}$ is fiber multiplication by a nonzero constant $c\in \O^*_\X(\X_i\cap
\X_j\cap\X_k)$.

The automorphism $\alpha_{ijk}$ naturally induces an automorphism $\beta_{ijk}$
of the restriction of $\N_i$ to $\X_i\cap
\X_j\cap\X_k$, where
$$
\beta_{ijk}=f_i^*\Bigl(\det (f_{i*}\alpha_{ijk})^{\otimes a-b}\otimes 
\det (f_{i*}\alpha_{ijk}\otimes id_{\omega _{f_i}})^{\otimes
b}\Bigr)
$$
and one easily checks that, by (\ref{rew}), $\beta_{ijk}$ is fiber multiplication by  $c^{-1}$.
We conclude that the line bundles $\poin _i\otimes \N_i\in \Pic \X_i$ can be glued together to
a line bundle $\poin$ over $\X$. It is clear that $\poin$ satisfies the thesis 
(since the $\poin _i$'s
do so).
\end{proof}
\begin{remark}
\label{nofine}
If the condition $(d-g+1,2g-2)=1$ is not satisfied the scheme $\pd$ can still be
constructed (as in the first part of the proof of \ref{fine}).
By \ref{cap} (\ref{capgeo})   $\pd$  is a geometric GIT-quotient if and only
if 
$(d-g+1,2g-2)=1$; if such a condition does not hold, there   exists  an open subset
$\mgb^d$ of $\mgb$ over which $\pd$ (and  $\pdg$) restricts to a geometric quotient. 
Such a nonempty open subset $\mgb^d$ is precisely the
locus of $d$ -general curves
by \ref{cap} (\ref{capstab}).
\end{remark}
\begin{nota}
\label{}

An application of  Lemma~\ref{poincare}  gives the existence of the analog of a
 Poincar\'e line bundle for the compactified Picard variety
of a family of automorphisim-free stable curves. More precisely, let $f:\X \la B$
be such a family and let
$\mu:B\la \mgb ^0$ be its moduli map; assume that $(d-g+1,2g-2)=1$. Then  we can form the 
compactified Picard scheme
$$
\pfb:=B\times_{\mgb ^0} \pdg \la B.
$$
Now, 
on the open subset of $\pdg$ lying over $\mgb^0$
there is a tautological curve 
$\taut$
which is constructed exactly as  ${\mathcal C}_{\pd}$ over $\pd$
(cf. proof of \ref{fine}). Observe that $\taut$ is a family of quasistable
(not stable) curves.
We can pull back $\taut$ to $\pfb$ and obtain a tautological curve
$\taut_f:=B\times_{\mgb ^0} \taut\la \pfb$.

Lemma~\ref{poincare} yields the analog of the Poincar\'e line bundle on $\taut$ and
hence  on
$\taut_f$; some care is needed  as the boundary points of $\pdg$ correspond to equivalence
classes of line bundles  that disregard the gluing data over the
exceptional component (see \ref{pdgbar} and \ref{glu} for the precise statement).

The construction of 
Poincar\'e line bundles over compactified Jacobians is an interesting problem in its own
right;
a  solution  within the category of
algebraic spaces was
 provided by E. Esteves in \cite{esteves}
applying  different techniques from ours. 

As we  indicated, our method
allows us to construct Poincar\'e  bundles  for automorphism-free  curves.
Rather than filling in the above missing details, 
we  
``stackify" the  construction of \cite{caporaso}  so that some of our results
will generalize to all stable curves (with or without automorphisms)
\end{nota}

\begin{nota}
\label{stack}
Let us introduce the stacks
defined by the group action used above:
$$\pdgst:=[H_d/G] \  \text{ and } \  \pdst:=[H_d^{st}/G]$$
When are they   Deligne-Mumford stacks (in the sense of of \cite{DM} and \cite{vistoli})? Do they
have a modular description? We begin with the first question,  adding to the picture the
 ``forgetful" morphisms
to $\mgbst$. To define it, pick  a scheme $T$ and
a section of $\pdgst$ (or of $\pdst$) over  $T$, that is a pair $(E\la T, \psi)$
where $E$ is a $G$-torsor and $\psi:E\la H_d $ is a $G$-equivariant morphism. Then we apply \ref{notpoin}
to obtain a family
$\X \la T$ of quasistable curves; 
the forgetful morphism maps $(E\la T, \psi)$ to the stable model of $\X \la T$
 (the reason why we call it ``forgetful" will be more clear from \ref{modstack}).

A map of stacks
${\mathcal P}\la {\mathcal M}$ is called {\it representable} (respectively, {\it strongly representable})
if given any algebraic space (respectively, scheme) $B$
with a map $B\la {\mathcal M}$, the fiber product $B\times _{\mathcal M} {\mathcal P}$ is
an algebraic space (respectively, a scheme).
\end{nota}

\begin{thm}
\label{DMstack}
 The stacks $\pdst$ and $\pdgst$ are
 Deligne-Mumford stacks if and only if
$(d-g+1,2g-2)=1$. In that case
the natural  morphisms
$\pdst \la \mgbst$ and  $\pdgst \la \mgbst$ are strongly representable.
\end{thm}
\begin{proof} 
As already said in \ref{cap} and in the proof of \ref{fine}, $H_d/G$ and $H_d^{st}/G$ are
geometric GIT-quotients 
(equivalently all stabilizers  are finite and reduced)
if and only if $(d-g+1,2g-2)=1$.
Hence the first sentence follows from the well known fact that a quotient stack 
like ours is a  Deligne-Mumford stack
if and only if all stabilizers are  finite and reduced.

For the second sentence, 
we first apply a common criterion for  representability 
(see for example \cite{AV}
4.4.3): our morphisms are representable if for every algebraically closed field $k'$ and every
section 
$\xi$ of $\pdst$ (respectively of $\pdgst$) over 
$\Spec k'$ the automorphism group of   $\xi$ injects into the automorphism group
of its image $X$ in $\mgbst$. This follows  from \ref{cap} (\ref{capstab}):
 $\xi$ is a map onto a $G$-orbit in $H_d$
and $\Aut(\xi)$ the stabilizer of such orbit (up to isomorphism, of course); the curve $X$ is
the stable model of the  projective curve $Z$
corresponding to such orbit, hence \ref{cap} (\ref{capstab}) gives us the desired injection.

We  obtained that the two forgetful morphisms   in the statements are representable, hence if $B$ is any
scheme and $B\la \mgbst$ the map corresponding to a family of curves $f:\X \la B$, 
the fiber product
$$
\pfb:=B\times _{\mgbst}\pdgst
$$
is an algebraic space;
  it remains to show that $\pfb$ is a scheme
(the fact that $B\times _{\mgbst}\pdst$ is also a scheme follows in the same way, or
observing that it is an open subspace of $\pfb$).
To do that, fix $\mu_f:B\la \mgb$  the moduli map of $f$ and consider the scheme
$$
Q_f:=B\times _{\mgb}\pdg
$$
which is projective over $B$ (if the fibers of $f$ are free from automorphisms
then $Q_f=\pfb$).
We shall  prove that there is a (natural) finite  projective
morphism
$$
\rho:\pfb\la Q_f;
$$
hence $\pfb$ is a scheme   (cf. \cite{viehweg} 9.4) projective over $B$.

To define $\rho$  we use
 \cite{vistoli} section 2 (in
particular 2.1 and 2.11),
which gives us that $\mgb$ and $\pdg$ are  the coarse moduli
schemes of $\mgbst$ and $\pdgst$ respectively  and that we have a canonical commutative diagram where
$\pi$ and $\pi'$ are proper
\begin{equation}
\begin{array}{lcccr}
 &\pdgst&\stackrel{\pi}{\la} &\pdg & \\
&\downarrow& & \downarrow\   & \\
B\la &\mgbst &\stackrel{\pi'}{\la} &\mgb&
\end{array}
\end{equation}
The two above maps from $B$ to $\mgbst$ and $\mgb$ are the same defining $\pfb$ and
$Q_f$; we let $\rho$ to be the base change over $B$ of   $\pi:\pdgst \la \pdg$, so that $\rho$ is
proper.

Now let $\lambda \in Q_f$ be a closed point.
Two different points in $\rho^{-1}(\lambda)$ correspond to two different maps
$
\psi, \psi':G\la H_d
$
mapping onto the orbit determined by $\lambda$,  hence (just as before) $
\psi$ and $ \psi'$ correspond to a nontrivial element in the stabilizer of a point in that
orbit. Since stabilizers are finite  $\rho$ has finite fibers; as $\rho$ is proper
we are done.
\end{proof}
\begin{nota}{\it Geometric description of $\pdst$ and $\pdgst$.}
\label{desc}
The modular description of $\pdst$ and $\pdgst$ can be  given by directly interpreting the
quotient stacks that define them; what we are going to obtain
 is a rigidified ``balanced Picard
stack".  The definition of the Picard scheme as a moduli scheme
representing a certain functor, or a certain stack, 
is well known to require   care, in fact a
subtle ``sheafification" procedure is needed to achieve representability.
The crux of the matter is that line bundles always possess automorphisms that fix the scheme
they live on, namely, fiber multiplication by nonzero constants; 
 see for
example
\cite{BLR} chapter 8 and
\cite{ACV} section 5. 
We are here in a  fortunate situation as the stacks  already exist and have 
some good properties (by \ref{DMstack}), all we have to do is to give them  a geometric
interpretation.

By \ref{cap} (\ref{caphigh}) we are free 
to assume  that $d$ is
very large.
 
Begin with an object in $\pdst$ (respectively in $\pdgst$), so
let $E\la T$ be a principal $PGL(r+1)$-bundle and $\psi:E\la H_d^{st}$
(respectively $\psi:E\la H_d$)
an equivariant map. Pulling back to $E$ the universal polarized family 
over the Hilbert
scheme
we obtain a  polarized family of stable (respectively quasistable) curves 
 over $E$, denoted as in \ref{notpoin} by
($f_E:\ZZ_E \la E , \L_E$). By construction $G=PGL(r+1)$ acts freely and we can form the quotient
$f:\X=\ZZ_E/G\la E/G=T$ which is a family of stable (respectively quasistable)
curves. 
Applying lemma~\ref{poincare} we obtain a balanced line bundle $\poin\in \Pic \X$ of relative
degree
$d$. Notice that $\poin$ is determined up to 
tensor product by pull-backs of line bundles on $T$, note also  that, using \ref{va}, we have a natural
isomorphism
$E\cong PGL(\pr {} (f_*\poin))$.

 Conversely let $(f:\X \la T, \L)$ be a pair consisting of  a family $f$ of stable (respectively
quasistable) curves and a balanced line bundle  of relative degree $d$ on $\X$; we now invert the previous
construction by
 producing  a principal $G$-bundle
$E\la T$ and a $G$-equivariant  map $E\la H_d^{st}$ (resp. $E\la H_d$).
 We 
 argue similarly to \cite{edidin} 3.2.
By \ref{va} $\L$ is relatively very ample and $f_*\L$ is locally free of rank $r+1=d-g+1$;
let $E\la T$ be the principal $PGL(r+1)$-bundle associated to 
the $\pr{r}$-bundle ${\mathbb P} (f_*\L)\la T$.
To obtain the equivariant map  to the Hilbert scheme consider the pull-back family
$f_E:\X_E=E\times _T\X\la E$ polarized by the balanced, 
relatively very ample line bundle $\L_E$ (pull-back of $\L$).
By construction ${\mathbb P} (f_{E*}\L_E)\cong \pr{r}\times E $
so that $\X_E$ is isomorphic over $E$ to a family of projective  curves in $\pr{r}\times E $
embedded by the balanced line bundle $\L_E$. By the universal property of the Hilbert
scheme this family determines a map
$\psi:E\la \Hilb$ whose image is all contained in $H_d^{st}$
(respectively in $H_d$). 
It is obvious that $\psi$ is $G$-equivariant.
\end{nota}
\begin{nota}
\label{modstack}
Let us summarize the construction of the previous paragraph, assume that
$(d-g+1,2g-2)=1$, then
\begin{enumerate}[{\bf (1)}]
\item
\label{balstack}
The stack $\pdst$
is the  ``rigidification" (in the sense of \cite{ACV} 5.1, see \ref{aut} below)
of the category whose sections over a scheme $T$
  are pairs
$(f:\X \la T, \L)$ where $f$ is a family  of stable curves of genus $g$ and $\L\in \Pic\X$
is  a balanced line bundle of relative degree $d$.
The arrows  between two such pairs are given by   cartesian diagrams
\begin{equation}
\begin{array}{lccr}
 & \X&\stackrel{h}{\la} &\X' \\
&\downarrow & &\downarrow\\
& T&\la &T'
\end{array}
\end{equation}
 and $\L\cong h^*\L'\otimes f^*M$ for $M\in \Pic T$.
\item
\label{qstack}
The stack $\pdgst$ is the  rigidification of the category whose sections over a scheme $T$ are
pairs
$(f:\X \la T, \L)$ where $f$ is a family  of quasistable curves of genus $g$ and $\L\in \Pic\X$
is  a balanced line bundle of relative degree $d$ .
Arrows are defined exactly as 
 in (\ref{balstack}).
\end{enumerate}

\begin{remark}
\label{aut} The rigidification procedure
removes those automorphisms of an $\L$ that fix $\X$; this is necessary  for
representability over $\mgbst$ (cf. \ref{DMstack} and \cite{AV}
4.4.3).
\end{remark}
In \cite{pandha} section 10, the scheme $\pdg$ was given a geometric 
 description in terms of rank-one torsion free sheaves rather  than line bundles.
This should enable one to obtain an alternative geometric description 
of the stacks $\pdst$, $\pdgst$ (and, obviously, of the scheme $\pd$).
\end{nota}
\begin{nota}
\label{fibre} Assume that $(d-g+1,2g-2)=1$ and 
let $f:\X \la B$ be a family of stable curves of genus $g$;
consider the schemes (cf. \ref{DMstack})
$$
\pf=B\times _{\mgbst}\pdst \   \text{ and }\  \pfb=B\times _{\mgbst}\pdgst .
$$
{\it If $(d-g+1,2g-2)\neq 1$
the two schemes $\pf$ and $\pfb$ can be defined in exactly the same way, provided that 
every singular fiber of $f$   is $d$-general.}

In fact, by \ref{cap} (\ref{capstab}),
the points in $H_d$ lying over the open subset $\mgb^d$ of $\mgb$,
parametrizing $d$-general curves, 
 are   all GIT-stable. Therefore
 the analogue of \ref{DMstack} holds, simply by restricting the quotient groupoids
over $\mgb^d$
(the proof is the same).

In the special case  $B=\Spec k$, so that the family $f$ reduces to a fixed stable curve
 $X$, we naturally change the notation and denote by $\pX$ (respectively by $\pXb$) 
the fiber of
$\pdst$ (respectively of $\pdgst$ ) over
$X$ as above.

 $\pX$ is a finite disjoint union of isomorphic copies of the generalized Jacobian of
$X$; the union is  parametrized by the set of
stably balanced multidegrees. 
Since   $X$  is $d$-general  a multidegree is balanced if and only if it
is stably balanced and  every $\delta \in \dcg^d$ has a unique balanced
representative (by \ref{dgen}).
Therefore
\begin{equation}
\label{fibreq}
\pX  \cong \coprod _{\md \in \domX^d}\picX{\md} \cong 
\coprod _{\delta\in \dcg^d}\picX{\delta}.
\end{equation}
The next result generalizes \ref{fine}.
\end{nota}

\begin{cor}
\label{rep}
Let $f: \X \la B$ be a family of 
 stable curves and    $d$ an integer. Assume that  
  every singular fiber of $f$ is $d$-general.
Then
the functor $\pffun$ is coarsely represented by the separated scheme
$\pf$;  if $B$ is regular, $\pf$ is smooth over $B$.
\end{cor}
\begin{remark}
\label{more} Under the  assumption that $(d-g+1,2g-2)=1$ the proof shows that $\pf$ is 
a fine
moduli scheme. 
\end{remark}
\begin{proof}  If we assume $(d-g+1,2g-2)=1$ the statement
follows   from \ref{DMstack} and \ref{modstack} 
and we 
obtain (as stated in \ref{more}) that $\pf$ is a fine moduli space. If, more generally,
 the singular fibers of
$f$ are $d$-general,
 we are still in the locus  where the
quotient defining $\pdst$ is geometric
(cf. \ref{fibre}). Then the statement follows as before
(the  reason why we get only a coarse moduli space is that
the   Poincar\'e line
bundle has been constructed only under the hypothesis
 $(d-g+1,2g-2)=1$). 
$\pf\la B$ has  equidimensional nonsingular fibers
(cf. (\ref{fibreq}) above), hence $\pf$ is smooth over $B$.
\end{proof}

\section{N\'eron models and balanced Picard schemes}
\label{comparison}
With the notation introduced in \ref{fibre},  we are ready to prove
our parametrization result.
\begin{thm}
\label{picner}
Let $f: \X\la B$ be a family  of stable curves 
of genus $g\geq 3$ such that $\X$ is regular
and $B$ is a one-dimensional regular connected scheme with function field $K$.
Let $d$ be such that every singular fiber of $f$ is $d$-general
(for example, assume that $(d-g+1,2g-2)=1$). 

\begin{enumerate}[(i)]
\item
\label{picner1}
Then $\pf$ is  the N\'eron model  of $\jacd$ over $B$.
\item
\label{picnersec}
If $f$ admits a section, 
 $\pf $ is isomorphic
to the N\'eron model \  $\ner$ of the Jacobian of the generic fiber of $f$.
\end{enumerate}
\end{thm}
\begin{proof}
If $f$ admits a section then $\Pic_K^d\cong \Pic ^0_K$ hence
$\NN(\Pic_K^d) \cong\NN( \Pic ^0_K)$. Thus the second part of the statement is an 
immediate consequence of the first.

By \ref{rep} $\pf$ is a smooth separated scheme of finite
type over $B$;  by \cite{BLR}  1.2/Proposition 4
it suffices, for part (\ref{picner1}),
to prove that $\pf$ is a {\it local}
N\'eron model, that is, we can replace $B$ by $\Spec R$ where $R$ is the
local ring 
 of $B$  at a
closed point (hence a discrete valuation ring of $K$).
Thus, we shall  assume that $f:\X \la \Spec R$ with $\X$
regular.
By \ref{sep} we have 
$$
\NN(\Pic_K^d)= \frac{\coprod _{\delta \in \dcg^d}\picf{\delta}}{\sim}
$$
(where ``$\sim$" denotes gluing along the generic fiber).

Since the closed fiber $X$ is $d$-general, a multidegree $\md$ is balanced if and only if it is
stably balanced and there is a natural bijection betweed the set 
of balanced multidegrees $\domX^d$
and $ \dcg^d$
 (cf. \ref{dgen}). Therefore we have a canonical $B$-isomorphism
$$\NN(\Pic_K^d)\cong
 \frac{\coprod_{\md \in \domX^d}\picf{\md}}{\sim}
$$
We now claim that there is canonical $B$-isomorphism 
\begin{equation}
\label{pfclaim}
\pf \cong \frac{\coprod_{\md \in \domX^d}\picf{\md} }{\sim}
\end{equation}
which, comparing the last two identities,  concludes the proof.

To prove (\ref{pfclaim}) it suffices to observe that the both schemes 
represent the balanced Picard 
functor for the given family $f$: for $\pf$ this follows from \ref{rep}, for the right hand
side this is clear
\end{proof}

\begin{remark}
In \ref{picner} the hypothesis that $\X$ is regular is necessary, see
 \ref{nosing} for an example illustrating why.
\end{remark}
We can apply the previous result to compare at least birationally different completions of
the generalized Jacobian. 
\begin{cor}  
\label{compare}
Under the same hypotheses of \ref{picner} (\ref{picnersec}),
let ${\overline \Pic_K^0}$ be any completion of $\Pic_K^0$ over $B$.
Then there exists a 
regular map (canonical for any fixed group structure on $\pf$) from the smooth
locus of 
${\overline \Pic_K^0}\la B$ to $\pf$, which restricts to an isomorphism on the generic fiber.
\end{cor}
\begin{proof}
Apply the N\'eron mapping property 
to $\pf$ (which we can do by \ref{picner}) and the unicity of the N\'eron model.
\end{proof}
\begin{remark}
\label{many}
It has been known for a long time that there is more than one good way of completing the
generalized Jacobian of a family of nodal (reducible) curves.
Perhaps the first to observe and study this phenomenon were T. Oda and C.S. Seshadri in \cite{OS};
their paper only dealt with a fixed curve and not with a family,
nevertheless 
 the  insights  contained there have deeply influenced the subsequent work of many authors.

Since then, diverse techniques have led to different models of compactified Jacobians.
The problem remains as to which completions are more suitable for
the miscellany of mathematical problems in which a compactified Picard variety is needed;
the previous result may be viewed in this perspective, 
offering a way of comparing different
constructions in different degrees.

A remarkable case is $d=g-1$,  which  has been particularly studied
(partly in relation with the problem of extending the theta-divisor).
Some correlation results have been proved by V. Alexeev in
\cite{alexeev} where there is also an overview of  the various
existing constructions.
As mentioned in \ref{g-1rm}, the  $d=g-1$ case is  ``degenerate" from our point of view
(arguing as \ref{degjac}, the compactified Picard variety is seen to have fewer components than
the N\'eron model).  For some aspects, however,
it turns out to be easier to handle precisely because of  certain degeneracy phenomena.

\begin{example}
\label{degjac} 
The previous corollary  applies  to the compactified Jacobians
given by the fibers of $\pdg$ over curves that that are not  $d$ -general.
For any family $f:\X \la B$ of (automorphism-free) stable curves of genus $g$  denote, as usual,
$
{\overline\pf}:=\pdg\times _{\mgb}B
$ and note that
 ${\overline\pf}$ depends on $d$, in fact the fibers of $\pdg$
over
$\mgb$ depend on $d$, as we are going to illustrate. If $X$ is a singular fiber of $f$,
the fiber of $\pdg$ over $X$ is denoted $\pXb$.

The simplest case in which we find a ``degenerate"
compactification of the generalized Jacobian is $d=0$
(this example works similarly if $d=g-1$). Let $X=C_1\cup C_2$ with $\#(C_1\cap C_2)=k$
and assume, which is crucial, that $k$ is even. 
Now,
$\dcg =\Z/k\Z$ and  the class
$$
\delta :=[(-\frac{k}{2},\frac{k}{2})]=[(\frac{k}{2},-\frac{k}{2})]
$$
has two balanced representatives (the ones above). Correspondingly, in 
${\overline{P_X^0}}\subset {\overline{P_{0,g}}}$,
line bundles having such multidegrees are strictly GIT-semistable and get identified 
to points having a stabilizer of positive dimension
 (the  so-called ``ladders", curves  obtained by blowing up all the nodes of $X$,
see \cite{caporaso} 7.3.3 for  details). Therefore the corresponding component of the
N\'eron model,
$\picX{\delta}$ (cf.\ref{sep}),  does not appear as an irreducible component of
${\overline{P_X^0}}$, where it collapses to a positive codimension boundary stratum.

In fact ${\overline{P_X^0}}$ has $k-1$ irreducible components, each of which corresponds to one of
the remaining classes in $\dcg$. Thus \ref{compare} implies that
if $f$ and $d$ are as in \ref{compare},  with $X$ as closed fiber, there is a diagram of
birational maps
\begin{equation}
\begin{array}{lccr}
f: & \overline{P^0_f}&\da &\overline{P^d_f} \\
&\uparrow & &\uparrow \\
&P^0_f&\ha& \pf
\end{array}
\end{equation}
\end{example}
\end{remark}
and the lower horizontal arrow is not an isomorphism.

\begin{nota}
\label{type}
Let $f:\X \la \Spec R$ be a family of generically smooth  curves
with closed fiber $X$ reduced, nodal and connected (not necessarily stable).
Let $\ner$ be the N\'eron model of its Jacobian; then its special fiber $\ner _k$
only depends on the geometry of $\X$,
or, which is the same, on the intersection form defined on the
minimal desingularization of $\X$
(see \cite{lorgroup}, \cite{edix} and \cite{BL} for explicit details and computations). More
precisely, the total space $\X$ can only have rational
singularities of type $A_n$ 
(i.e. formally equivalent to $xy=u^{n+1}$)
at the nodes of $X$, and the  singularities that
will interfere with the structure of $\ner _k$ are those occurring at the {\it external} nodes of $X$
(i.e.  nodes lying on two different components).
Let $\delta$ be the number of external nodes of $X$ and suppose that
$\X$ has a singularity of type  $A_{n_i}$ at the $i$-th external node.
Then the structure of
$\ner _k$ only depends on $\underline{n}=(n_1,\ldots, n_{\delta})$ so that we can denote
$
\NX^{\underline{n}}
$ the special fiber of a N\'eron model of this type.

 We need  the case where $\X$ is nonsingular, so that 
$\underline{n}=(0,\ldots, 0)$; then we  denote the special fiber of the N\'eron model
of the Jacobian of $f$ by
$$
\NX:= \NX^{(0,\ldots, 0)}
$$
We have for any nodal (connected) curve $X$ (see \ref{sep})
\begin{equation}
\label{NX}
\NX\cong \coprod _{\delta \in \dcg}\picX{\delta}
\end{equation}
\end{nota}

\begin{example}
\label{nosing}
We  now exhibit an example showing that the assumption that $\X$ be regular in \ref{picner}
cannot be weakened  by assuming $\X$ normal.
Let $f:\X\la \Spec R$ having as closed fiber  $X=C_1\cup C_2$ with 
$k=\#(C_1\cap C_2)\geq 2$.
Assume that $\X$ has a singularity of type $A_n$ 
 at one of the nodes of $X$
and it is smooth otherwise. Then the twister group $\twX$ of $f$
is generated by $T_1:=\O_{\X}((n+1)C_1)\otimes \O_X$
which has multidegree $\mdeg T_1=(-(nk+k-n), nk+k-n)$. Thus the group of multidegree
classes for such an $f$ will be
(using a notation similar to the one introduced in \ref{type})
$$
\dcg^{(n,0,\ldots,0)}\cong \Z/(nk+k-n)\Z
$$
which is bigger than $\dcg$ (if $n\geq 1$ of course).
The closed fiber 
$
\NX^{(n,0,\ldots,0)}
$ of the
N\'eron model of the generalized jacobian of $f$ has component group  isomorphic to 
$\Z/(nk+k-n)\Z$, whereas the components of the closed fiber of $\pf$ are parametrized
by $\dcg$ (if $X$ is $d$- general). 

Finally, if $X$ is not $d$- general
so that we are in a degenerate case as described in  \ref{degjac}, the
number of components of the special fiber of $\pf$ is  smaller than $\#\dcg$ and hence also
smaller than  $\#\dcg^{(n,0,\ldots,0)}$.
\end{example}

\begin{nota}
\label{}

A natural side question is: when are $\pdg$ and $\overline{P}_{d',g}$ isomorphic?
Similar question for  the  stacks.
This is easy to answer, we do it  for the schemes
but it is obvious that the same answer holds for the stacks. By  \ref{cap} (\ref{caphigh}) we have that 
$\pdg\cong \overline{P}_{d',g}$
if and only if $d\pm d'\equiv 0 \mod (2g-2)$ and these isomorphisms are canonical.
Then we just need to count; denoting ``$\Phi$"  the Euler $\phi$-function on natural numbers we
have
\end{nota}
\begin{lemma}
\label{euler}
The number 
 of non isomorphic $\pdg$ for which $(d-g+1,2g-2)=1$ is equal to
$\Phi (g-1)$ if $g$ is odd and
to $\frac{\Phi (g-1)}{2}$ 
if $g$ is even
\end{lemma}
\begin{proof}
As we said,
 there are exactly $g$ non isomorphic models for  $\pdg$.
We  choose as representatives for each class of such models the values for $d$ given by
$
d=0,1,\ldots,g-1
$
so that we have
$$
\overline{P}_{0,g}\cong \overline{P}_{2g-2,g},\  \  \overline{P}_{1,g}\cong \overline{P}_{2g-3,g},\ldots,
\overline{P}_{g-2,g}\cong \overline{P}_{g,g}
$$
and for any $d'\geq 2g-2$  
$$\overline{P}_{d',g}\cong \overline{P}_{-d',g}\cong \overline{P}_{e,g}$$
where $0\leq e<2g-2$ and
$d'=n(2g-2)+e$.
Now
 $(d-g+1,2g-2)=1$ implies  $(d,g-1)=1$; 
if $g$ is odd, one immediately sees that the converse  holds, 
and we are done.

If $g$ is even,  the condition $(d-g+1,2g-2)=1$
is equivalent to $d$  even  and  coprime with $g-1$. So the  values of $d$  that we are
counting are the positive even integers
$d$ coprime with $g-1$  and smaller than $g-1$.
This number  equals $\frac{\Phi (g-1)}{2}$ (just notice that for any odd 
$m\in {\mathbb N}$, the Euler function $\phi(m)$ counts an equal number of odd and even integers; 
in
fact if 
$r$ is
odd and coprime with $m$, the even number $m-r$ is also coprime with $m$; 
same thing starting with
$r$ even.)
\end{proof}

\section{Completing N\'eron models  via N\'eron models}
\label{completion}

\begin{nota}
\label{coprime}
 From now on we shall assume that the stable curve
$X$ is $d$-general (\ref{dgen}).
  For example, one may assume that $(d-g+1,2g-2)=1$.
\end{nota}

Fix  $f:\X \la B=\Spec R$  a family of stable curves with  smooth
 generic fiber and regular total space $\X$. In \ref{fibre}   we introduced the scheme 
$\pfb$,
projective over
$B$ which, by \ref{picner},  is  a compactification of the N\'eron model of the
Picard variety $\Pic ^d_K$ (by \ref{picner}); recall that 
$\pXb$ denotes its closed fiber.
In the present section we shall exhibit a stratification
of  $\pXb$ in terms of 
   N\'eron models associated to all the
 connected partial normalizations of $X$ (Theorem \ref{strata}).
In section \ref{quotient} we shall prove that $\pfb$ is dominated
 by the  N\'eron model of 
a degree-$2$ base change
of 
$\jacd$.
See 
\cite{andreatta} for a different approach to
the problem of compactifying  N\'eron models of Jacobians.

\begin{nota}
\label{pdgbar}
With the notation  introduced in \ref{fibre},
we shall refer to the points in $\pXb \smallsetminus \pX$ 
as the ``boundary points of $\pXb$".
To describe them precisely  we  need some simple preliminaries.

Let $X$ be a stable curve, the quasistable curves of $X$ (cf. \ref{notstab})
correspond
bijectively to  sets of its nodes: let $S$ be a set of nodes of $X$,  we shall denote
$
\nu _S:\XS \la X
$
the normalization of $X$ at the nodes  in $S$ and 
$$
\YS := \XS \cup (\bigcup _1^{\#S}E_i)
$$
the quasistable curve of $X$ obtained by joining the 
two points of $\XS$ lying over the $i$-th
node in $S$ with a smooth rational curve
  $E_i\cong
\pr{1}$ (so that one may call $\YS$ the {\it blow up of $X$ at $S$}).
\end{nota}
\begin{nota}
\label{glu}
A point  of $\pXb$ corresponds to an equivalence class
of  pairs $(\YS,L)$ where $S\subset
X_{sing}$
 and $L\in \Pic^d \YS$ is a balanced line bundle.
Two pairs $(\YS,L)$ and $(\YS',L')$ are equivalent if and only if
$\YS = \YS '$ and $L_{|\XS}\cong L'_{|\XS}$.

The boundary points are those for which $S\neq\emptyset$

\end{nota}

\begin{remark} 
\label{connected}
Notice that  a quasistable curve $\YS$ of $X$
admits a (stably) balanced line bundle
(of degree $d$) if and only if the subcurve $\XS$ (obtained by removing all of the exceptional
components) is connected.

In fact if $\XS=Z_1\cup Z_2$ with $Z_1\cap Z_2=\emptyset$ then a stably balanced $\md$
has to satisfy
$
d_{Z_1\cup Z_2}=m_{Z_1\cup Z_2}
$,
on the other hand
$d_{Z_1\cup Z_2}=d_{Z_1}+d_{Z_2}$ and hence $d_{Z_1}=m_{Z_1}$
(and $d_{Z_2}=m_{Z_2}$). This is impossible  as the complementary curve 
of $Z_1$, containing $Z_2$,
is not a union of exceptional components (cf \ref{BIrm}).
\end{remark}

\begin{nota}
\label{one}
Fix the quasistable curve $\YS$ and consider $\dcY$; recall that a balanced multidegree
must have degree $1$ on all exceptional components of $\YS$, so that
not all elements in $\dcY$ have a balanced representative.
Denote 
$$
\dcYY :=\{\delta\in \dcY : \  \delta \  \text{has a balanced representative }\}
$$ 
Thus for every
 $\delta \in \dcYY$ there exists a unique
(by \ref{coprime}) balanced representative
which we shall denote

\begin{equation}
\label{breq}
(d_1^{\delta},\ldots,d_{\gamma}^{\delta},1,\ldots,1)
\end{equation}
so that $[(d_1^{\delta},\ldots,d_{\gamma}^{\delta},1,\ldots,1)] =\delta$ and 
$\sum _1^{\gamma}d_i^{\delta}=d-s$, where  $s:=\#S$.

By \ref{connected}
we have that $\dcYY $ is empty if and only if $\XS$ is not connected.

\end{nota}
The next lemma will be used in the proof of Theorem~\ref{strata}.
\begin{lemma}
\label{key}
Using the above notation, assume $\XS$ connected. Then the map

$$
\rho : \dcYY\la \dcS,  \hskip.5in 
[(d_1^{\delta},\ldots,d_{\gamma}^{\delta},1,\ldots,1)] \mapsto
[(d_1^{\delta},\ldots,d_{\gamma}^{\delta})]
$$
is bijective.
\end{lemma}
\begin{proof}
As we said $\rho$ is well defined because of  the assumption \ref{coprime}.
We shall use the notation of \ref{domlm} and \ref{domb},
together with the following:
let $Z\subset \XS\subset Y$,
set
$k_Z^S:=\#(Z\cap {\overline{\XS\smallsetminus Z}})$
and denote by
$e_Z$
the number of points in which $Z$ meets the the exceptional components
of $\YS$
so that 
\begin{equation}
\label{kez}
k_Z=e_Z+k_Z^S.
\end{equation}
The map $\rho$ can be factored as follows:

\begin{equation}
\label{keyeq}
\begin{array}{lccccccr}
\rho : &\dcYY&\la &\dom_{\YS}(\mb_{\YS}^d)&
\stackrel{\pi}{\la} & \dom_{\XS}(\mb) &\stackrel{\sigma}{\la} &\dcS \  \   \\
\\
&\delta & \mapsto & (d_1^{\delta},\ldots,d_{\gamma}^{\delta},1,\ldots,1)&\mapsto
&(d_1^{\delta},\ldots,d_{\gamma}^{\delta})&\mapsto
&[(d_1^{\delta},\ldots,d_{\gamma}^{\delta})]
\end{array}
\end{equation}
where $\mb =(b_1,\ldots, b_{\gamma})$ with
$$
b_i:=\frac{d}{2g-2}w_{C_i}- \frac{e_{C_i}}{2}
$$
and  $w_{C_i}=2g_{C_i}-2+k_{C_i}$.

To prove that $\rho$ is surjective, first of all observe that,
by \ref{domlm}, $\sigma$ is surjective.
Now we claim that given $\md =(d_1,\ldots, d_{\gamma},1,\ldots,1)\in \Z^{\gamma +s}$ 
such that
$|\md|=d$, we have that 
$\md$ is balanced if and only if
for every $Z\subset \XS$ we have 
\begin{equation}
\label{Zext}
m_Z(d)\leq d_Z\leq M_Z(d)-e_Z
\end{equation}
where $M_Z(d)= \frac{d}{2g-2}w_{Z}+ \frac{k_Z}{2}$ and $m_Z(d)=M_Z(d)-k_Z$ as usual.
In fact for every exceptional component $E$ of $\YS$ 
we have 
$
w_Z=w_{E\cup Z}
$
and hence the basic inequality on $Z\cup E$ gives
\begin{displaymath}
d_Z+1=d_{Z\cup E}\leq\left\{ \begin{array}{ll}
M_Z(d)+1 &\text{ if } (E\cdot Z)=0\\
M_Z(d) &\text{if } (E\cdot Z)=1\\
M_Z(d)-1 &\text{if } (E\cdot Z)=2\\
\end{array}\right.
\end{displaymath}
Iterating for all $E$ we get the claim.

Therefore
$\md$ is balanced if and only if (using (\ref{kez}))
$$
\frac{d}{2g-2}w_{Z}- \frac{k_Z^S}{2}-\frac{e_Z}{2}\leq d_Z\leq \frac{d}{2g-2}w_{Z}+
\frac{k_Z^S}{2}-\frac{e_Z}{2}
$$
if and only if
$$
(d_1,\ldots, d_{\gamma})\in \dom_{\XS}(\mb )
$$
This shows that $\rho$ is surjective; to prove that it is injective
it suffices to show that $\sigma $ is (the other two arrows of diagram (\ref{keyeq})  are obviously
injective). If $ \dom_{\XS}(\mb) $ contains two equivalent multidegrees,
then, using  (\ref{Zext}), we would get that there exists a subcurve $Z\subsetneq \XS$
for which $m_Z(d)$ is integer, which is impossible
(as usual, by assumption \ref{coprime}).
\end{proof}

\begin{nota}
\label{pfb}
By \ref{DMstack} and \ref{modstack},\  
$\pfb$ is a coarse moduli scheme for the functor from $B$-schemes to sets
which associates to a $B$-scheme $T$ the set of equivalence classes of
pairs $(h:\Y\la T, \L)$
where $h:\Y\la T$ is a family of quasistable curves having $\X _T$ as stable model; and $\L$ is
a balanced line bundle on $\Y$. The equivalence relation is the same as in \ref{balfun}.
\end{nota}
\begin{nota}
\label{nod}
The structure of the closed fiber $\pXb$ of $\pfb$
does not depend on $d$ (by \ref{coprime})
and  is a good compactification of $\NX$ (see \ref{type}). Therefore
we shall introduce the notation 
$$\NXb:=\pXb$$
Such a completion can be described 
by means of the N\'eron models of the Jacobians of all connected partial normalizations of $X$:
\end{nota}
\begin{thm}
\label{strata} $\NXb$ has a natural stratification as follows
\begin{equation}
\label{strataeq}
\NXb \cong \coprod_{\stackrel{S\subset X_{sing}:}{\XS \text{connected}}} \NS
\end{equation}
Denote $Q_S\subset \NXb$ the stratum  isomorphic to $\NS$ under the  decomposition
(\ref{strataeq}); then
\begin{enumerate}[(i)]
\item
\label{stratacod}
$Q_S$ has pure codimension 
$\#S$.
\item
\label{strataord}
$Q_S\subset \overline{Q_{S'}}$ if and only if $S'\subset S$.
\item
\label{stratasm}
The smooth locus of \  $\NXb$ is $\NX$.
\end{enumerate}
\end{thm}
\begin{proof}
As we explained in \ref{pdgbar}, the points  of $\pXb=\NXb$ parametrize  pairs $(\YS,L)$ in such a way that
for every $S\subset X_{sing}$ we have a well defined locus $Q_S$  in $\pXb$,  corresponding to 
balanced line bundles  on $\YS$.
For example, $\pX$ corresponds to the stratum $S=\emptyset$ (isomorphic to
$\NX$).

 In turn, $Q_S$ is a disjoint union of irreducible components 
isomorphic to the generalized Jacobian of $\XS$
(cf. \ref{glu} and \ref{fibre}); there is one component for every (stably) balanced
multidegree on
$\YS$. More precisely, 
for any balanced $\md=(d_1,\ldots, d_{\gamma},1,\ldots,1)$ on $\YS$ let us denote $\md^S=(d_1,\ldots,
d_{\gamma})$  its restriction to $\XS$. 
Then the moduli morphism
$$
\Pic^{\md}\YS \la \pXb
$$
(associated to the universal line bundle on $\Pic^{\md}\YS\times \YS$, see \ref{pfb})
factors through a surjective morphism followed by a canonical embedding
\begin{equation}
\label{stratafact}
\Pic^{\md}\YS \stackrel{}{\twoheadrightarrow} \Pic^{\md^S}\XS\ha Q_S\subset \pXb
\end{equation}
(see \ref{glu})
whose image is open and closed in $Q_S$.

Set $\delta^S:=[\md^S]\in \dcS$.
We shall now see that the components of $Q_S$
are in one-to-one correspondence with the elements of $\dcS$.
The balanced multidegrees on $\YS$ are bijectively parametized by $\dcYY$
(cf. \ref{one});  by \ref{key}
the restriction to $\XS$ of a balanced multidegree induces the bijection
$$\rho:\dcYY\leftrightarrow \dcS$$ of \ref{key},
so we  are done.
In other words we obtain the stratification in the statement of our Theorem
\begin{equation}
\label{!!}
Q_S   \cong
\coprod _{\delta^S \in \dcS}\Pic ^{\delta^S}\XS \cong \NS
\end{equation}
where the second isomorphism is (\ref{NX}).

Part (\ref{stratacod}) is a simple dimension count. We already know that each irreducible component
of
$Q_S$ is isomorphic to the generalized Jacobian of $\XS$; the genus of $\XS$ is equal to $g-s$
 hence we are done.

By the previous results, part(\ref{strataord}) 
follows from Proposition 5.1 of \cite{caporaso} (see below for more details).

Now  (\ref{stratasm}); 
 quite generally, the N\'eron mapping property applied to \'etale points
implies that  any completion $\overline{N}$ of a N\'eron model $N$ over $B$
must be singular along $\overline{N}\smallsetminus N$
(If $\overline{N}\smallsetminus N$ contained regular points one would use  2.2/14
of
\cite{BLR} and find an \'etale point of $N_K$ which does not come from an \'etale
point of $N$).
We include a direct proof to better illustrate the structure of $\NXb$.

It suffices to prove that every component of every positive codimension
stratum is contained in the closure of more than one irreducible component of $\NX=\pX$.
This  also follows from  Proposition 5.1 of \cite{caporaso}. Let us treat the case $\#S=1$;
then $\YS$ has only one exceptional component $E$ intersecting (say)
$C_1$ and $C_2$ (viewed now as components of $\XS$ by a slight abuse of notation).
If the point $(\YS, L)$ belongs to the component of $Q_S$ corresponding to the multidegree
$(d_1,d_2,\ldots, d_{\gamma},1)$,  we have that
$(\YS, L)$ is contained in the closure of the two components of $\pX$ that correspond to
multidegrees $(d_1+1,d_2,\ldots, d_{\gamma})$ and $(d_1,d_2+1,\ldots, d_{\gamma})$.
\end{proof}
\begin{nota}
\label{Gx}
Let $X$ be a stable curve; as we have seen, $\NXb$ has
a stratification 
(by equidimensional, possibly disconnected  strata) parametrized by the sets of nodes of $X$
which do not disconnect
$X$, denote by $\glu$ this set:
$$
\glu:=\{S\subset X_{sing}: \XS \  \text{is connected} \}
$$

For some more details on the stratification of Theorem~\ref{strata},
introduce the dual graph $\G$  of $X$, (cf. \ref{comb})
and recall the genus formula
$
g=\sum _1^{\gamma}g_i + b_1(\G)
$
where $g_i$ denotes the geometric genus of $C_i$ and $b_1(\G)$ is the first Betti number
(see \ref{betti}).
\end{nota}
\begin{cor} 
\label{min}Let 
$X$ be a stable curve and $S\in \glu$; let
$Q_S\subset \NXb$ be a stratum as defined in Theorem~\ref{strata}. Then
\begin{enumerate}[(i)]
\item
\label{mindim} $\dim Q_S\geq \sum _1^{\gamma}g_i$
\item
\label{min1}
$\dim Q_S =\sum _1^{\gamma}g_i \Longleftrightarrow \XS$ is of compact type $\Longleftrightarrow Q_S$
is irreducible.
\item
\label{minnum}
The number of minimal strata of $\NXb$ (in the stratification of Theorem~\ref{strata})
is equal to  $\#\dcg$.
\end{enumerate}
\end{cor}
\begin{proof}
(\ref{mindim})  is equivalent to 
$
\dim Q_S \geq g-b_1(\G),
$
hence, by \ref{strata} (\ref{stratacod}),
 it suffices to show that $\#S\leq b_1(\G)$. Thus we must prove  that the maximum number of nodes of
$X$ that can be normalized without disconnecting
the curve is $b_1(\G)$.
Equivalently,   that the maximum number of edges of
$\G$ that can be removed without disconnecting
 $\G$ is $b_1(\G)$.
This follows from \ref{betti}.

Now we prove (\ref{min1}).
$\dim Q_S =\sum _1^{\gamma}g_i $ if and only if $Q_S$ is a minimal stratum of $\NXb$ (by
\ref{strata} and part (\ref{mindim})), if and only if all the nodes of $\XS$ are separating (i.e. any
partial normalization of $\XS$ fails to be connected), if and only if
$\XS$ is if compact type (by definition, cf.
\ref{ct}). This proves the first double arrow of part (\ref{min1}).

$\XS$ is if compact type if and only if its dual graph is a tree, if and only if
$\Delta_{\XS}=\{0\}$ (this can be easily shown directly or it follows from \ref{degcomp}),
if and only if $Q_S$ has  only one irreducible component
(by \ref{strata}\  $Q_S\cong \NS $ whose components correspond to elements in $\Delta_{\XS}$).
This concludes (\ref{min1}).

Now (\ref{minnum}). The  strata of minimal dimension (equal to $\sum _1^{\gamma}g_i $)
correspond bijectively to
the connected partial normalizations of $X$ that are of compact type which,
in turn, correspond (naturally) to the spanning trees of the dual graph of $X$
(cf. \ref{st}).
Now, the number of spanning trees of  $\G$
(the so called ``complexity" of the graph) is shown to be equal to
the cardinality of $\dcg$ in \ref{degcomp}. So we are done.
\end{proof}
\begin{example}
\label{vine}
Let $X=C_1\cup C_2$ with $C_i$ nonsingular and $\#(C_1\cap C_2)=k$;
then the set
$\glu$ is easy to describe:
$\glu =\{S\subset X_{sing}: S\neq X_{sing}\}$.
Given $S\in \glu$ let $\#S=s$ so that
$\XS=C_1\cup C_2$ with $\#(C_1\cap C_2)=k-s$.

The  connected 
components of $\NX$,
each isomorphic to the generalized jacobian of $X$,
 are parametrized by $\Z/k\Z$.

The strata $Q_S$ of codimension $1$ of $\NXb$ are parametrized by the nodes of $X$,
denoted $n_1,\ldots,n_k$. If $S=\{n_i\}$, $Q_{n_i}$ is the 
special fiber $\NS $ of the N\'eron model of the Jacobian
of a family specializing to the normalization of $X$ at $n_i$;  hence it is made of $k-1$
connected  components of dimension $g-1$.

And so on, going down in dimension till the minimal strata,
which correspond to the $k$ curves of compact type obtained from $X$ by normalizing
it at $k-1$ nodes.
Each of these strata is isomorphic to the 
closed fiber of the
N\'eron model of the Jacobian of a specialization to a curve of compact type
having $C_1$ and $C_2$ as irreducible components; therefore
it is an irreducible projective variety (isomorphic to $\Pic^{\underline 0}C_1\times 
\Pic^{\underline 0}C_2$) of dimension $g-k+1$.
\end{example}

\section{The compactification as a quotient}
\label{quotient}

We begin with some informal remarks to motivate the content of this last section;
consider a family of nodal curves $f:\X \la B=\Spec R$ having regular $\X$ and
singular closed fiber $X$. 
Let $p\in X$ be a nonsingular point, then $p$ corresponds to a degree-$1$ line
bundle of $X$ which, up to an \'etale base change of $f$
(ensuring the existence of a section through $p$) is the specialization of a degree
$1$ line bundle on the generic fiber. So $p$ corresponds to a unique point
in $\NN(\Pic _K^1)$.

What if $p$ is a singular point of $X$? Of course (intuitively)
$p$ can still be
viewed as a limiting configuration of line bundles on $X$.  On the other hand there
will never be a section
passing through $p$ (not even after \'etale base change of $f$). What is needed to
have such a section is a ramified base change, in fact a degree-$2$ base change
will suffice 
(because $X$ has ordinary double points).
If $f_1:\X _1\la B_1$ is the base change of $f$ under a degree-$2$ ramified covering
$B_1=\Spec R_1\la B$, then $\X _1$ has a singularity of type $A_1$ at each node of
the closed fiber $X_1\cong X$. If $p_1\in X_1$ is the point corresponding to $p$,
then $f_1$ (or some \'etale base change) does admit a section through $p_1$ ,
therefore $p_1$, and hence our original point $p$, corresponds to a unique point of
$\NN(\Pic ^1_{K_1})$.

All of this suggests that to complete the N\'eron model of the Picard variety of
$\X_K$ we could use the N\'eron model of the the Picard variety of a ramified base
change of order $2$.
To better handle the N\'eron models $\NN(\Pic ^d_{K_1})$ we shall introduce and
study the minimal desingularization of $\X_1$, 
whose closed fiber is the quasistable
curve
$Y$ of $X$ obtained by blowing up all the nodes of $X$.

\begin{nota}
\label{ladder}
Let $X$ be a stable curve; 
consider  
the quasistable curve $Y$
obtained by blowing up all the nodes of $X$ so that,
with the notation of \ref{pdgbar}, $Y:=Y_{X_{sing}}$.
Denote
$$
\sigma: Y \la X
$$
the morphism contracting all of the exceptional components of $Y$.

Recall now that, by \ref{strata},  $\NXb$ has
a stratification labeled by  $\glu$.
We shall  exhibit a   decomposition
 of $\NY$  labeled by $\glu$
and  prove that it is  naturally related to
the stratification of $\NXb$.

By \ref{dgen}  for any $\delta \in \dcYd$ there exists a unique semibalanced
representative
$\md^{\delta}$.
Fix now a set $S$ of nodes of $X$ and define
$$
\dcYS:=\{\delta \in \dcYd: d^{\delta}_E =1  \Leftrightarrow \sigma(E)\in S\}
$$
\end{nota}

Let $\gamma$ be the number of irreducible components of $X$ and let $s=\#S$;
order the
exceptional components of $Y$ so that the first $s$ are those corresponding to $S$
(i.e. mapped to $S$ by $\sigma$).
Connecting with   \ref{one}  we can partition the component group $\Delta_Y$
of $\NY$  using $\glu$:

\begin{lemma}
\label{keyY} Let $Y=Y_{X_{sing}}$.
\begin{enumerate}[(i)]
\item
\label{}
For every $S$ there is a natural bijection 
\begin{equation}
\label{bij}
 \  \dcYS \leftrightarrow \dcYY, \hskip.2in 
[(d_1^{\delta},\ldots,d_{\gamma}^{\delta},1,\ldots,1,0,\ldots,0)] \mapsto
[(d_1^{\delta},\ldots,d_{\gamma}^{\delta},1,\ldots,1)]
\end{equation}
\item
\label{keydec}
$
\coprod _{S\in \glu} \dcYS= \dcYd
$
\end{enumerate}
\end{lemma}
\begin{proof}
Let $\md =(d_1,\ldots,d_{\gamma},1,\ldots,1)$ a multidegree on $\YS$ and denote
its ``pull-back" to $Y$ by 
$\md^*=(d_1,\ldots,d_{\gamma},1,\ldots,1,0,\ldots,0)$; to prove that (\ref{bij}) is a
bijection it suffices to prove that
 $\md$ satisfies the basic inequality on $\YS$ if and only if
$\md^*$ satisfies the basic inequality on $Y$.
Denote $\sigma_S:Y\la \YS$ the contraction of all exceptional components of $Y$ that do not
correspond to $S$. Let $Z\subset Y$ be a subcurve and denote
$Z_S=\sigma_S(Z)\subset \YS$.
Then it is easy to see that
$w_Z=w_{Z_S}$ and that
$k_{Z}=k_{Z_S}+2t_{Z}$ where $t_Z$ is the number of exceptional components $E$
of $Y$ that are
not contained in $Z$ and such that $\#(E\cap Z)=2$.
If we write the basic inequality for $Z_S\subset \YS$ as usual
(omitting the dependence on $d$ which is fixed)
\begin{equation}
\label{ZS}
m_{Z_S}\leq d_{Z_S}\leq M_{Z_S}
\end{equation}
the basic inequality for $Z\subset Y$ is
\begin{equation}
\label{dd}
m_{Z_S}-t_Z\leq d^*_Z\leq M_{Z_S}+t_Z.
\end{equation}
Under the correspondence (\ref{bij}) we have  $d_{Z_S}=d^*_Z$, hence
 it is obvious that, if $\md$ satisfies (\ref{ZS}), then $\md^*$ satisfies (\ref{dd}).
Conversely, suppose that $\md^*$ satisfies the basic inequality and let 
$Z_S\subset \YS$ be a
subcurve. Denote by $Z=\sigma ^{-1}(Z_S)$ so that 
$t_Z=0$; thus  the basic inequality for $Z_S$
is  the same as for $Z_S$ and hence $\md$ satisfies it.

For the second part, recall that, because of \ref{connected},
  $\dcYY$ is empty if and only
if
$S\not\in
\glu$ (see \ref{one}). Thus $\dcYS$ is empty if $S\not\in \glu$ and the second part of the
lemma follows.
\end{proof}
\begin{remark}
\label{}
As a consequence we get the  $\glu$-decomposition  of $\NY$   mentioned in \ref{ladder}: 
$$
\NY = \coprod_{S\in \glu}\bigr(\coprod_{\delta \in \dcYS}\Pic ^{\delta}Y\bigl)
$$
\end{remark}

\begin{nota}
\label{two}
Let $f:\X \la \Spec R=B$ with $\X$ regular and assume that $f$ admits a section.
The curve $Y$ 
(defined in \ref{ladder})
is the closed fiber of the regular minimal model of the base
change of
$\gen$ under a degree-$2$ ramified covering of $\Spec R$. More precisely, let $t$ be a
uniformizing parameter of $R$ and let
$K\ha K_1$ be the degree-$2$ extension $K_1=K(u)$ with $u^2=t$. Denote $R_1$ 
the DVR of  $K_1$ lying over $R$,
so that $R\ha R_1$ is a degree $2$ ramified extension.
Denote
$B_1=\Spec R_1$ and consider the  covering $B_1\la B$.
The corresponding base change of $f$ is denoted 
$$
f_1:\X_1:=\X \times _BB_1\la B_1
$$
and $X_1$  its closed fiber. At each of the nodes of $X_1$ the total space $\X _1$ has a
singularity formally equivalent to 
$
xy=u^2,
$
which can be resolved by blowing up once each of the nodes of $X_1$ (see \cite{DM} proof of
1.2). Denote $\Y \la \X_1$ this blow-up and $h:\Y \la B_1$ the composition;
thus $h$ is a family of quasistable curves having $\X_1$ as stable model and $Y$ as closed fiber.
We summarize with a diagram
\begin{equation}
\begin{array}{lccr}
\;  \;    \Y&\la &\X \\
{\tiny h}\downarrow  & &\   \downarrow {\tiny f}\\
\;  \;   B_1 &\la &B
\end{array}
\end{equation}
Denote  $\pich{d}\la B_1$ the Picard variety for $h$ and $\jacdh$ its generic fiber.
\end{nota}
\begin{prop}
\label{quot} In  the set up of \ref{two},
let $\NN(\jacdh)\la B_1$ be the N\'eron model of $\jacdh$; then there is a canonical
surjective
$B$-morphism
$$
\pi: \NN(\jacdh)\la \pfb.
$$
The restriction of $\pi$ to the closed fibers is compatible with their
$\glu$-stratifications in the following sense: for any $S\in \glu$ the
restriction $\pi _S$ of $\pi$
is a surjective   morphism 
$$
\pi _S: \coprod_{\delta \in \dcYS}\Pic ^{\delta}Y\la Q_S\cong
\NS
$$
(notation of
\ref{strata}) all of whose closed fibers are isomorphic to $(k^*)^s$ with $s=\#S$.
\end{prop}

\begin{remark}
\label{}
$\pi$ is described as a quotient by a torus action in
\ref{torus}.
\end{remark}
\begin{proof}
By \ref{sep} we have
$
\NN(\jacdh)\cong \frac{\coprod _{\delta \in \dcYd}\pich{\delta}}{\sim}.
$
The crux of the proof is to show that for every $\delta \in \dcYd$
 there is a canonical
morphism
$$
\mu_{\delta}:\pich{\delta}\la \pfb .
$$
To do that, let $S$ be the unique element in $\glu$ such that
$\delta \in \dcYS$ and consider the unique semibalanced representative 
$\md ^{\delta}$ of $\delta$ 
(cf. \ref{dgen}).  Denote by $T$  and identify (by \ref{identify})
$
T:=\pich{\delta}=\pich{\md^{\delta}}.
$
Set
$$
h_T:\Y_T=\Y\times _{B_1}T\la T
$$
and let $\Poi$  be the Poincar\'e line bundle on $\Y_T$.
Now we  apply the construction of \ref{polmod} to $h_T=p$ and $\Poi=\N$. 
Thereby we obtain a family,
$
\overline{\Y_T}\la T
$
 (by contracting all the exceptional components of the fibers
of $h _T$ where $\Poi$ has degree equal to zero) and a  line bundle $\L$  on 
$\overline{\Y_T}$ which pulls back to $\Poi$.
The singular closed fibers of ${\overline{\Y_T}}\la T$
are all isomorphic to $\YS$ and $\L$ has balanced multidegree $\md =
(d_1^{\delta},\ldots,d_{\gamma}^{\delta},1,\ldots,1)$ (the fact that $\L$ is balanced
follows from the proof of \ref{keyY}, whose notation we are here using). It may be useful to
sum up  the construction in a diagram where all squares are cartesian:
\begin{equation}
\begin{array}{lccccccr}
\overline{\Y_T}  & \longleftarrow &\Y_T&\la &\Y &&\\
\downarrow & &\downarrow & &\downarrow&&&\\
\X_T&=&\X_T&\la&\X_1&\la& \X\\
\downarrow &&\downarrow & &\downarrow & &\downarrow\\
T&=&\pich{\delta}&\la&B_1&\la &B
\end{array}
\end{equation}
Now the pair $(\overline{\Y_T}\la T, \L)$ is a family of quasistable curves with a
balanced line bundle of degree $d$. The stable model of $\overline{\Y_T}$ is
$\X_T$ therefore (by \ref{pfb}) we obtain a moduli morphism
$$
\mu _{\delta}:T=\pich{\delta} \la \pfb.
$$
As $\delta$ varies, the morphisms $\mu _{\delta}$ agree on the smooth fibers,
that is, away from the closed point of $B$. Therefore
(as in the proof of \ref{picner}) they glue together
to a $B$-morphism $ \pi: \NN(\jacdh)\la \pfb$
as stated.

To prove the rest of the statement it suffices to look at the
closed fiber, as $\pi _K$ is obviously a surjection,
in fact
$$
\NN(\jacdh)_{K_1}= \jacd \times _B\Spec K_1\ =(\pfb)_K\times _B\Spec K_1=(\pf)_K\times
_B\Spec K_1
$$
Now by \ref{keyY} and \ref{key} 
(and with the same notation) we have natural bijections
\begin{equation}
\begin{array}{lccccr}
 \dcYS&\leftrightarrow&\dcYY&\leftrightarrow &\dcS \  \  \  \  \\
\\
\md^{\delta}&\mapsto &\md = [(d_1^{\delta},\ldots,d_{\gamma}^{\delta},1,\ldots,1)]
&\mapsto &\md^S = [(d_1^{\delta},\ldots,d_{\gamma}^{\delta})]
\end{array}
\end{equation}
As we said, the  singular fibers of $\overline{\Y_T}\la T$ are isomorphic
to $\YS$ and we  proved above that the restriction of $\mu_{\delta} $ to
the closed fibers factors
$$
\Pic^{\md^{\delta}}Y \stackrel{\cong}{\la}
\Pic^{\md}\YS\stackrel{}{\twoheadrightarrow}
\Pic^{\md^S}\XS\ha  \pXb
$$
where we  used (\ref{stratafact}) for the last two arrows;
the rest of the proof naturally continues as  that of  \ref{strata}.
\end{proof}

\begin{nota}
\label{torus} Let $b=b_1(\G)$.
It is not difficult at this point to interpret $\pi$ as a quotient by a natural
 action
of $(k^*)^b$ on $\NY$ (extended to a trivial action on $\NN(\jacdh)$).
Observe that $\Pic ^{\md}Y\cong \Pic ^{\md}\YS \cong \Pic^{\md^S}X $
(notation in the proof of \ref{strata})
and that $b-s=b_1(\Gamma _{\XS})$;
denote $X^{\nu}$ the normalization of $X$,  
we have a diagram of canonical exact sequences
\begin{equation}
\begin{array}{lccccccccr}
 &&0& &0& &\\
 &&\downarrow & &\downarrow & &\\
0&\la & (k^*)^s & =&(k^*)^s&\la &0 &&\\
& &\downarrow & &\downarrow&&\downarrow &\\
0&\la & (k^*)^b&\la &\Pic^{\md}\YS&\stackrel{\nu^*}{\la}&
\Pic^{\md^S}X^{\nu}\times \prod_1^s\Pic^1\pr{1}&\la& 0\\
 &&\downarrow & &\downarrow & &\downarrow\\
0&\la & (k^*)^{b-s}&\la &\Pic^{\md^S}\XS&\stackrel{\nu^*}{\la}&
\Pic^{\md^S}X^{\nu}&\la& 0\\
 &&\downarrow & &\downarrow & &\downarrow\\
 &&0 & &0 & &0\\
\end{array}
\end{equation}
(where $\nu^*$ always denotes pull-back via the normalization map).
The middle vertical sequence describes  the restriction of
$\pi_S$ to any irreducible component, $\Pic^{\md}\YS$,
as the quotient of the action of $(k^*)^s $
on the gluing data over the exceptional components of $\YS$.

\end{nota}

We applied  the following standard fact
(included  for completeness).
\begin{lemma}
\label{polmod}
Let $p:\ZZ \la T$ be a family of semistable curves 
of genus at least $2$ over a scheme $T$. Let $\N\in \Pic
\ZZ$ having non-negative degree on each exceptional component of the fibers of $p$.
Then there exist
\begin{enumerate}[(a)]
\item
 a factorization of $p$ 
$$
p:\ZZ \stackrel{\psi}{\la}\overline{\ZZ}\stackrel{\overline{p}}{\la} T
$$
via a family of semistable
curves $\overline{p}$ and a birational morphism $\psi$ which contracts some exceptional
components of the fibers of $p$;
\item
a line bundle $\NL\in \Pic \overline{\ZZ}$ having positive degree on all exceptional
components of the fibers of $\overline{p}$ and such that
$\psi ^*\NL\cong \N\otimes p^*M$, where $M\in \Pic T$. 
\end{enumerate}
\end{lemma}
\begin{proof}

For $n$ high enough (how high depends on $\N$) we have that
$\omega _p^n\otimes \N$ is relatively base-point-free and
$p_*(\omega _p^n\otimes \N)$ is a vector bundle on $T$
(trivial variation on  Corollary to Theorem 1.2 in \cite{DM} p.78).
Moreover $\omega _p^n\otimes \N$ defines a birational morphism  
$\psi:\ZZ \la \overline{\ZZ}\subset {\mathbb P} (p_*(\omega _p^n\otimes \N))$ 
contracting the exceptional components of $p$ where $\N$ has degree $0$.
The line bundle $\NL$ is given by $\NL = \O _{\overline{\ZZ}}(1)\otimes \omega
_{\overline{p}}^{-n}$.
\end{proof}
\begin{remark}
\label{reldet}
It is clear that 
$\overline{\ZZ}$ is uniquely determined (just contract all the exceptional components of the fibers of
$p$ where
$\N$ has degree $0$) whereas $\overline{\N}$
is determined
only up to pull-backs of line bundles on $T$. More precisely the lemma gives a map
form $\Pic \ZZ/p^*\Pic T \la \Pic\overline{\ZZ}/{\overline p}^*\Pic T$.
\end{remark}

\begin{remark}
\label{}
We conclude by observing  that, as a consequence of \ref{two},
 the completion 
$\pfb$ of the N\'eron model satisfies a
mapping property for smooth schemes
defined over quadratic, possibly ramified,  coverings
of $B$. This should be viewed as a strengthening of
the mapping property of N\'eron models
with respect to smooth schemes 
defined over \'etale coverings of $B$.
It is in fact well known (see \cite{artin}
section 1)  that  N\'eron models are functorial with
respect to
\'etale base changes, but not in general.

To be more precise, let $Z$ be a scheme smooth over $\Spec R_1$
(where $R\ha R_1$ is a ramified quadratic
extension as in \ref{two}), and let $v_K:Z_{K_1}\la \jacd$
be a $K$-morphism. Then there exists a unique
$B$-morphism $v:Z\la \pfb$ extending $v_K$. Of course $v$ is obtained
by first extending the lifting  of $v_K$ to $Z_{K_1}$
$$u_{K_1}: Z_{K_1}\la \NN(\jacdh)_{K_1}= \jacd \times _B\Spec K_1,$$
by the N\'eron mapping property
$u_{K_1}$ extends to $u:Z\la \NN(\jacdh)$; thus $v$ is the composition
of $u$ with
$\pi$ (defined in \ref{quot}).
\end{remark}

\section{Appendix}
\label{abel} 
This appendix is made of two distinct parts.
The first   illustrates  some  applications of the
results in the paper.
The second  part 
 summarizes some well known  combinatorial facts which have been used throughout.

\

\centerline{\it Applications: towards Brill-Noether theory of stable curves}

\begin{nota}
\label{sembal}
Let $f:\X\la B$ be a  family of stable curves and 
$T$  a scheme over $B$.
Let $p:\ZZ \la T$  be a family of semistable curves having $\X_T$ as stable
model;
if $\L\in \Pic \ZZ$ is balanced
of relative degree $d$,
 we can associate to $\L$ a unique map
$$
\mu_{\L}:T \la \pfb, \hskip.3in t\mapsto [\L _{|p^{-1}(t)}]
$$
which we call the {\it moduli map} of $\L$ (note that in this case $\ZZ$ is
necessarily  a family of quasistable
curves).

More generally, suppose that $\N\in \Pic \ZZ$ is semibalanced (cf. \ref{BI}). Apply the construction of 
\ref{polmod} to  obtain a pair
$(\overline{\ZZ}, \NL)$ (so that 
$\overline{\ZZ}\la T$ has  $\X_T$ as stable model).
Then $\NL\in \Pic \overline{\ZZ}$ is a balanced line bundle
and its moduli map $\mu_{\NL}:T \la \pfb$ can be viewed as induced by 
$\N$.
In summary, to a semibalanced line bundle on $\ZZ$ one associates a unique map
$T\la \pfb$.
\end{nota}

\begin{nota}
\label{AN}
Let $f:\X \la \Spec R=B$ be a family of curves with $\X$ regular
and reducible closed fiber $X$ (as usual),
denote $f_d:\X^d_B\la B$ its $d$-th fibered power.
Consider the degree-$d$ Abel map of the generic fiber
$$
\alpha_K^d:\gen^d \la \Pic^d\gen , \hskip.3in   (p_1,\ldots,p_d)\mapsto [\O_{\gen}(\sum p_i)];
$$
what is the limit of such a map as $\gen$ specializes to $X$?

   Not much is known about  defining (and completing)
  Abel maps for reducible curves.
A geometric construction   for irreducible  curves has  been carried out  in \cite{EGK}
building upon previous well known work of A.Altman and S.Kleiman.
Yet serious difficulties arise when reducible fibers occur
(even when restricting, as
we are, to nodal singularities).

As a first step towards understanding
Abel maps of reducible curves, we
consider  the  unique extension of $\alpha_K^d$
given by the N\'eron mapping property
$$
\alpha_f^d:\dX^d_B\la \NN(\Pic^d\gen)
$$
where $\dX^d_B=\X\smallsetminus sing(f_d)$; 
We refer to $\alpha _f^d$ as the
degree-$d$ {\it Abel-N\'eron  map} of $f$. The case $d=1$  has been studied by B.
Edixhoven   in
\cite{edix}, where there is also a characterization of when  
it is a closed immersion
(in the example  below it is).

The results of our paper enable us, on the one hand, to 
give a geometric description of the  Abel-N\'eron map
by identifying $\NN(\Pic^d\gen)\cong \pf$. On the other hand
we have a natural ambient space where one can construct a completion for it,
namely the compactification $\pfb$.

\end{nota}

\begin{example}
\label{ANex}
Fix a stable curve  $X=C_1\cup C_2$ with $C_1$ and $C_2$  smooth of  genus 
equal to $h\geq 1$ 
 and $\#C_1\cap C_2=2$ (thus
$g=2h+1$);
 let $f:\X \la \Spec R=B$ be a family of curves with $\X$ regular and $X$ as closed fiber.
Since   our $X$ is general for $1$,  we 
can identify $\NN(\Pic^1\gen)=P^1_f$ by (\ref{picner}) so that
the first Abel-N\'eron map becomes
$$
\alpha_f:\dX\la P^1_f.
$$
We claim   that

\begin{enumerate}[\bf (1)]
\item
\label{}
$\alpha_f$ can be completed to a map
$
\overline{\alpha_f}:\X\la \overline{P_f^1}
$;
\item
\label{} 
$\overline{\alpha_f}$ has a geometric description as the moduli map of a natural line bundle;

\item
\label{}
 the restriction to $X$ of $\overline{\alpha_f}$ 
 does not depend on $f$.
\end{enumerate}
Consider the line bundle $\L=\O_{\dX \times_B \X}(\Delta)\in \Pic (\dX \times_B \X)$
where $\Delta\subset \dX \times_B \X$ is the  diagonal in $\X^2_B=\X \times _B  \X$.
Then, applying the set up of \ref{sembal} with $T=\dX$, we  claim that
$\alpha _f$ is the moduli map of $\L$ (this is obviously true on the generic fiber $\gen$
of $\dX$). For that it suffices to show that $\L$ is balanced, i.e. that for every
nonsingular point $x\in X$ the line bundle
$\O_X(x)$ (the restriction of $\O_{\dX \times_B \X}(\Delta)$ to the fiber over $x$)
is balanced. This follows easily,   by checking  that for every subcurve
$Z$ of our $X$, we have $m_Z(1)<0$ and
$M_Z(1)>1$ so that we have
$$
\NN(\Pic^1\gen)\cong P^1_f\cong \frac{\Pic_f^{(0,1)}\coprod \Pic_f^{(1,0)}}{\sim}.
$$
 Now let us denote $r:\ZZ \la \X^2_B$ the resolution of singularities. A direct computation
shows that $\ZZ$ is obtained by replacing each of the four singular points of $\X^2_B$ by a
$\pr{1}$ so that
$p:\ZZ\la\X=T$ is a family of quasistable curves; 
moreover the proper transform $\tilde{\Delta}\subset
\ZZ$ of
 $\Delta$ defines a line bundle 
$\N=\O_\ZZ(\tilde{\Delta})$
having non-negative degree on  every exceptional component
of the fibers of $p$.
One checks that $\N$ is semibalanced
hence,
applying the construction of \ref{sembal},  we obtain a regular   map
$$\mu_{\NL}:T=\X \la \overline{P_f^1}$$ 
which defines the extension $\overline{\alpha _f}=\mu_{\NL}$ of $\alpha _f$ that 
we wanted.

To show that the restriction  $\alpha_X:=\overline{\alpha _f}_{|X}$
 does not depend on $f$
one simply observes that if $x\in X$ is a nonsingular point, then
its image is just the class of $\O_X(x)$. If $x$ is singular,
denote by $Y_x$ the quasistable curve obtained by blowing-up $X$ at $x$
and let $q\in Y_x$ be any nonsingular point of $Y_x$ lying in the unique exceptional component. 
Then, as $q$ varies, the line bundles $\O_{Y_x}(q)$  are all identified to the same
point $\lambda _x$ in $\overline{P_X^1}$
(by \ref{glu}); then the image 
$\alpha_X(x)$ is exactly the point  $\lambda _x$. 

We mention (without proof)
that  $\alpha_X$ is a closed 
embedding of $X$ into $\overline{P_X^1}$.
\end{example}
\begin{nota}
\label{}
The method of the previous example can be applied  to  all stable curves, but nontrivial
complications arise. First of all, it is not always true
 that the ``diagonal" line
bundle used above is balanced; a more delicate construction is needed to prove
that the same properties (1)-(3) hold.

The global version of such a morphism (mapping the universal curve over $\mgb$ to
$\overline{P_{1,g}}$) could also be carried out, as it is reasonable to expect, in view of the
independence on $f$ of the Abel-N\'eron map 
(property (3)). 

Let us finish with a few words about  the  Abel-N\'eron maps for higher degree $d$. 
The problem 
can be approached similarly to what outlined  for $d=1$;  however
the situation is  considerably more subtle.
One important difference is that, as soon as $d\geq 2$,
the  $d$-th Abel-N\'eron map will depend on $f$, for some combinatorially determined
cases. In other words, 
the analogue of
property  (3) fails.

Another difficulty is the fact (observed by E. Esteves) that a completion of the Abel map
will not be defined on $\X^d_B$, but only on
some modification $\widetilde{\X^d_B}\la \X^d_B$ of it.

 These
hurdles are  to be expected, as the set up
 leads towards a construction of Brill-Noether varieties for singular  curves.
As a first step, we can define the Brill-Noether scheme $\overline{W^0_d}(X,f)$ 
(generalizing the Brill-Noether variety  of effective  line bundles 
of degree $d$ on a smooth curve)
as follows: 
$$\overline{W^0_d}(X,f):=\overline{\im (\alpha_f^d)_{k}} \subset \pXb$$
 i.e. the closure of the image of the restriction
$
(\alpha_f^d)_{k}: \dot{X}^d\la \pXb 
$, where $\dot{X}$ denotes the smooth locus of $X$.
The closure symbol is used   because such a scheme  parametrizes 
``boundary points", that is, line bundles on quasistable curves $Y\neq X$ having $X$ as stable
model; we shall denote $W^0_d(X,f)$ its open subset parametrizing line bundles on $X$.

The presence of $f$  in the notation
is needed for $d\geq 2$;  although
 we can prove that $\overline{W^0_1}(X,f)$ never depends on
$f$, 
for  $d\geq 2$   this turns out to fail. 
To be more precise,
denote by $\X^{\nu}_{sep}$ the partial normalization of $X$ at 
 its separating
nodes (so that $\X^{\nu}_{sep}=X$ if $X$ has no separating node), then we 
conjecture the following. 
The restricted Abel-N\'eron map $(\alpha_f^d)_{k}$ is independent of $f$
if and only if every connected component of    $\X^{\nu}_{sep}$ is $k$-connected 
(i.e. admits no subset of $k$ disconnecting nodes) for every $k\leq d$.
\end{nota}

\centerline{\it Combinatorics of stable curves}

\begin{nota}
\label{comb}
Some   features of  stable curves are
nicely expressed using graph theory. 
Chapter 1 of the  article \cite{OS} contains a thorough 
study of combinatorial aspects of the theory of compactified Jacobians and
of degenerations of Abelian varieties. In the sequel we recall 
only a small number of facts that can be found in that paper.

To a nodal curve $X$ 
having $\gamma$ irreducible components and $\delta $ nodes,
one attaches a graph $\G$ defined as
the  symplicial complex (of dimension at most $1$)
defined to have  one {\it vertex}  for every irreducible component
of $C$, and one {\it edge}   connecting two vertices for every node in
which the two corresponding components intersect. 
Thus $\G$ has 
$\gamma$ vertices, $\delta$ edges and among the edges there is
a loop for every node  lying on a single irreducible component of $X$.
\end{nota}
\begin{nota}
\label{betti}
The first Betti number $b_1(\G)$ 
(sometimes called the {\it cyclomatic number})
is, for any orientation on $\G$
$$
b_1(\G):= \dim _{\Z}H_1(\G, \Z) =\delta -\gamma  +1
$$
Recall also 
that {\it the first Betti number of a connected graph is the maximal number of one-dimensional
open symplices that can be removed from the graph without disconnecting it}.

Another important, somewhat less standard,  invariant of a graph is its {\it complexity}

\begin{defi}
\label{st}
Let $\Gamma$ be a connected graph. A {\it spanning tree of $\G$} is a subgraph $\Gamma '\subset \Gamma$
which is a connected tree and such that $\Gamma $ and $\Gamma '$ have the same vertices.
The {\it complexity} of $\Gamma $,
$\comp (\Gamma)$, is defined to be the number of spanning trees that it contains.
\end{defi}
\end{nota}
\begin{example}
\label{ct} Let $X$ be connected.
\item(1) $X$ is of compact type if and only if $\G$ is a tree, if and only if
$b_1(\G) = 0$, if and only if $\comp (\G)=1$.

\item (2)
By the genus formula $g=\sum g_i + b_1(\G) $ we get that $b_1(\G) \leq g$.
Moreover, $b_1(\G) =g$ if and only if all irreducible components of $X$ have geometric genus $0$.
\end{example}

\begin{nota}
\label{hom}

The complexity can be computed cohomologically. Fix an orientation on $\Gamma$ and consider the standard
homology operators
\begin{equation}
\partial :  C_1(\Gamma , \Z) \la \ C_0(\Gamma , \Z), \hskip.3in
e \mapsto v-w
\end{equation}
where $e$ is an  edge of $\Gamma$,  starting in the vertex $v$ and ending in
the vertex $w$.
And 
\begin{equation}
\delta :  C_0(\Gamma , \Z) \la  C_1(\Gamma , \Z), \hskip.3in
v \mapsto \sum e^+_v -\sum e^-_v
\end{equation}
where  $e^+_v$ are the edges starting
at the vertex
$v$ and $e^-_v$ are those ending in $v$.
Then introduce the {\it complexity group} of the graph $\g$
$$
{\partial C_1(\Gamma , \Z)\over \partial \delta C_0(\Gamma , \Z)}
$$
the name ``complexity group" is due to the  theorem of Kirchhoff-Trent (\cite{OS} p.21)
stating that {\it such a group is finite and its cardinality is equal to the complexity of $\Gamma$.}

The next lemma 
is  Proposition 14.3 in \cite{OS}
(see also \cite{lorgraph}).
\end{nota}

\begin{lemma}
\label{degcomp} For a nodal connected curve $X$ with dual graph $\Gamma _X$ we have
$$
\dcg \cong {\partial C_1(\Gamma _X, \Z)\over \partial \delta C_0(\Gamma _X, \Z)}.
$$
In particular the cardinality of $\dcg$ is equal to the complexity of $\Gamma _X$.
\end{lemma}

\noindent Lucia  Caporaso \  - \  caporaso@mat.uniroma3.it\  \\\
 Dipartimento di Matematica, Universit\`a Roma Tre\\
Largo S.\ L.\ Murialdo 1 \  \  00146 Roma - Italy\\


\begin{thebibliography}{EGKH02}

\bibitem[AV01]{AV}
 D. Abramovich, A. Vistoli.
\newblock {\em Compactifying  the 
 space of stable maps.}
\newblock  Journal of the American Mathematical Society, 15 27-75, 2001.

\bibitem[ACV01]{ACV}
 D. Abramovich, A. Corti, A. Vistoli.
\newblock {\em Twisted bundles and admissible covers.}
\newblock  Comm. Algebra 31  Special issue in honor of Steven L. Kleiman.   (2003), 
no.8, 3547-3618


\bibitem[Al04]{alexeev} V. Alexeev: {\it Compactified Jacobians and Torelli map.}
Publ. RIMS, Kyoto Univ. 40 (2004), 1241-1265. 



\bibitem[An99]{andreatta} F. Andreatta: {\it On Mumford uniformization and N\'eron models of
Jacobians of semistable curves over complete rings.}   Moduli of abelian varieties (Texel
Island, 1999),  11-126, Progr. Math. 195, Birkh\"auser  Basel  2001


\bibitem[A86]{artin} M. Artin: {\it N\'eron models.}
Arithmetic geometry, edited by G. Cornell, J.H. Silverman. Proc. Storrs. Springer (1986)


\bibitem [BLR]{BLR}S. Bosch, W L\"uktebohmert, M. Raynaud:
\newblock N\'eron models.
\newblock Ergebnisse der Mathematik 21 Springer 1980.

\bibitem [BL02]{BL}
S. Bosch, D. Lorenzini:  {\it Grothendieck's pairing on component groups of
jacobians.}  Invent. Math. 148 (2002), 353-396.

\bibitem[C94]{caporaso}
 L. Caporaso.
\newblock {\em A compactification of the universal {P}icard variety over the moduli
 space of stable curves.}
\newblock  Journal of the American Mathematical Society, 7:589-660, 1994.

\bibitem[CCC04]{CCC}
 L. Caporaso, C. Casagrande, M. Cornalba:
\newblock {\em Moduli of roots of line bundles on curves.}
\newblock  To appear  in Transactions of the American Mathematical Society. MathAG/0404078


\bibitem[DM69]{DM}
P.Deligne, D.Mumford: 
\newblock {\it The irreducibility of the space of curves of given genus.}
\newblock Inst. Hautes \'Etudes
Sci. Publ. Math. No 36 (1969) 75-120.

\bibitem[E00]{edidin}
D.Edidin:
\newblock {\it  Notes on the construction of the moduli space of curves.}
\newblock in Recent Progress in Intersection Theori 
 Ellingsrud, Fulton,
Vistoli editors, Birkhauser 2000.


\bibitem[E98]{edix}
B.Edixhoven:
\newblock {\it  On N\'eron models divisors and modular curves.}
\newblock J.Ramanujan Math. Soc. 13 (1998), no2 157-194.

\bibitem[EGK00]{EGK}
E.Esteves, M. Gagn{\'e},  S. Kleiman.
\newblock {\em Abel maps and presentation schemes.}
\newblock Communications in Algebra 28(12) 5961-5992 (2000)

\bibitem[E01]{esteves}
E.Esteves.
\newblock {\em Compactifying the relative Jacobian over families of
   reduced curves.}
   \newblock  Transactions of the American Mathematical Society, 353(8):3045--3095, 2001.


\bibitem[F04]{fontanari}
C. Fontanari:
\newblock {\em On the geometry of the compactification of the universal
Picard variety}.
\newblock Rendiconti dell'Accademia Nazionale dei Lincei Roma (to appear)


\bibitem[Gie82]{gieseker}
 D.Gieseker:
\newblock  Lectures on moduli of curves.
\newblock  TIFR Lecture Notes 69 (1982) Springer.

\bibitem[SGA]{SGA} A. Grothendieck:
\newblock {\it Technique de descente et th\'eor\`emes d'existence en 
g\'eometrie alg\'ebrique. Le sch\'emas de Picard.}
\newblock S\'eminaire Bourbaki SGA 232, SGA 236


\bibitem[HM98]{HM}
 J. Harris, I. Morrison:
\newblock  Moduli of curves.
\newblock  Graduate texts in Math. 187 (1998) Springer

\bibitem[LM]{LM}
 G. Laumon, L. Moret-Bailly:
\newblock  Champs alg\'ebriques.
\newblock    Ergebnisse der Mathematik 39 Springer 2000.


\bibitem[L89]{lorgraph}
D. Lorenzini: {\it Arithmetical graphs.}  Math. Ann.  285 (1989), 481-501.

\bibitem[L90]{lorgroup}
D. Lorenzini: {\it Groups of components of N\'eron models of Jacobians.}
Comp. Math. 73 (1990), 145-160.

\bibitem[L93]{lorner}
D. Lorenzini: {\it On the group of components of a N\'eron model.}
J. reine angew. Math 445 (1993) 109-160.

\bibitem[M78]{maruyama}
M. Maruyama: {\it Moduli of stable sheaves II.}
J. Math Kyoto Univ. 18 (1978).


\bibitem[MR85]{MR}
N. Mestrano, S. Ramanan: {\it Poincar\'e bundles for families of curves.}
Journ. reine angew. Math. 362 (1985) 169-178.

\bibitem[M66]{mumford}
 D. Mumford:
\newblock  Lectures on curves of an algebraic surface.
\newblock  Annals of mathematics
studies Princeton University press 1966.



\bibitem[GIT]{GIT}
 D. Mumford,  J.  Forgarty, F. Kirwan:
\newblock  Geometric Invariant Theory.
\newblock   Third edition (First edition 1965)  Ergebnisse der Mathematik 34 Springer 1994.

\bibitem[N64]{neron}
A. N\'eron:
{\it Mod\`eles minimaux des vari\'et\'es ab\'eliennes sur les corps locaux et
globaux.}  Inst. Hautes \'Etudes Sci. Publ.Math. No. 21 1964 .

\bibitem[OS79]{OS}
T.Oda, C.S.Seshadri: {\it Compactifications of the generalized Jacobian variety.}
Trans. A.M.S. 253 (1979) 1-90

\bibitem[P96]{pandha}
R. Pandharipande.
\newblock {\em  A compactification over $\overline{M_g}$ of the universal moduli
  space of slope-semistable vector bundles.}
\newblock Journal of the American Mathematical Society 9(2):425-471,
  1996.

\bibitem[R70]{raynaud}
 M. Raynaud: {\it Sp\'ecialisation du foncteur de Picard.}  Inst. Hautes \'Etudes
Sci. Publ. Math. No. 38 1970 27-76.





\bibitem[S94]{simpson}
 C. T.   Simpson.
\newblock {\em Moduli of representations of the fundamental group of a smooth projective variety. }
\newblock  Inst. Hautes \'Etudes Sci. Publ. Math., 80: 5-79 , 1994.  



\bibitem[Vie91]{viehweg}
E. Viehweg:
\newblock  Quasi-projective moduli for polarized manifolds.
\newblock Ergebnisse der Mathematik (30) Springer.

\bibitem[Vi89]{vistoli}
A. Vistoli:
\newblock  {\em Intersection theory on algebraic stacks and their moduli spaces.}
\newblock Invent. math. 97 613-670 (1989).


\end{thebibliography}
\end{document}